\documentclass[11pt]{amsart}
\usepackage{amssymb,amsmath,a4wide}
\usepackage[pagebackref=true]{hyperref}
\usepackage[all]{xy}

\title{Towards a homotopy theory of higher dimensional transition
  systems}

\author[P. Gaucher]{Philippe Gaucher}

\address{CNRS UMR 7126\\Laboratoire PPS\\
  Univ Paris Diderot\\
  Sorbonne Paris Cit\'e\\ F-75205 Paris\\ France}

\urladdr{http://www.pps.jussieu.fr/{\~{}}gaucher/} 



\swapnumbers

\newcommand{\C}{\mathcal{C}}
\newcommand{\D}{\mathcal{D}}
\newcommand{\K}{\mathcal{K}}
\newcommand{\LL}{\mathcal{A}}
\newcommand{\I}{\mathcal{I}}
\newcommand{\W}{\mathcal{W}}
\newcommand{\F}{\mathcal{F}}

\newcommand{\de}{\partial}
\newcommand{\p}\times

\newtheorem{thm}{Theorem}[section]
\newtheorem{prop}[thm]{Proposition}
\newtheorem{lem}[thm]{Lemma}
\newtheorem{cor}[thm]{Corollary}

\newtheorem{defn}[thm]{Definition}
\newtheorem{nota}[thm]{Notation}
\newtheorem{rem}[thm]{Remark}
\newcommand{\bd}{\begin{defn}}
\newcommand{\ed}{\end{defn}}
\newcommand{\bp}{\begin{prop}}
\newcommand{\ep}{\end{prop}}
\newcommand{\bth}{\begin{thm}}
\renewcommand{\eth}{\end{thm}}
\newcommand{\bpf}{\begin{proof}}
\newcommand{\epf}{\end{proof}}

\newcommand{\fL}[1]{\ar@{->}[ll]_-{#1}}
\newcommand{\fR}[1]{\ar@{->}[rr]^-{#1}}
\newcommand{\fRr}[1]{\ar@{->}[rrr]^-{#1}}
\newcommand{\fD}[1]{\ar@{->}[dd]_-{#1}}
\newcommand{\fU}[1]{\ar@{->}[uu]^-{#1}}
\newcommand{\f}[2]{\ar@{->}[#1]|{#2}}
\newcommand{\ff}[2]{\ar@2{->}[#1]|{#2}}
\newcommand{\frr}[1]{\ar@{->}[rrrr]^-{#1}}

\newcommand{\fl}[1]{\ar@{->}[l]_-{#1}}
\newcommand{\fr}[1]{\ar@{->}[r]^-{#1}}
\newcommand{\fd}[1]{\ar@{->}[d]_-{#1}}
\newcommand{\fu}[1]{\ar@{->}[u]^-{#1}}

\newcommand{\iso}{\cong}

\renewcommand{\leq}{\leqslant}
\renewcommand{\geq}{\geqslant}
\newcommand{\dd}[1]{\uparrow\!\!{#1}\!\!\uparrow}

\def\cartesien{%
  \ar@{-}[]+R+<6pt,-2pt>;[]+RD+<6pt,-6pt>%
  \ar@{-}[]+D+<2pt,-6pt>;[]+RD+<6pt,-6pt>%
}
\def\cocartesien{%
  \ar@{-}[]+L+<-6pt,+2pt>;[]+LU+<-6pt,+6pt>%
  \ar@{-}[]+U+<-2pt,+6pt>;[]+LU+<-6pt,+6pt>%
}
\def\hocartesien{%
  \ar@{-}[]+R+<6pt,-2pt>;[]+RD+<6pt,-6pt>_{h}%
  \ar@{-}[]+D+<2pt,-6pt>;[]+RD+<6pt,-6pt>%
}
\def\hococartesien{%
  \ar@{-}[]+L+<-6pt,+2pt>;[]+LU+<-6pt,+6pt>_{h}%
  \ar@{-}[]+U+<-2pt,+6pt>;[]+LU+<-6pt,+6pt>%
}

\newcommand{\brm}[1]{\rm{\mathbf{#1}}}

\newcommand{\set}{{\brm{Set}}}

\DeclareMathOperator{\id}{Id}

\DeclareMathOperator{\Mor}{Mor}
\DeclareMathOperator{\pr}{pr}

\newcommand{\liminj}{\varinjlim}

\DeclareMathOperator{\whdts}{\brm{WHDTS}}

\DeclareMathOperator{\cts}{\brm{CTS}}

\DeclareMathOperator{\cub}{\underline{Cub}}

\DeclareMathOperator{\dom}{dom}
\DeclareMathOperator{\bl}{\brm{\underline{L}}}

\makeatletter
\def\varholim@#1#2{%
  \vtop{\m@th\ialign{##\cr
    \hfil$#1\operator@font holim$\hfil\cr
    \noalign{\nointerlineskip\kern1.5\ex@}#2\cr
    \noalign{\nointerlineskip\kern-\ex@}\cr}}%
}
\def\holimproj{%
  \mathop{\mathpalette\varholim@{\leftarrowfill@\textstyle}}\nmlimits@
}
\def\holiminj{%
  \mathop{\mathpalette\varholim@{\rightarrowfill@\textstyle}}\nmlimits@
}
\makeatother

\DeclareMathOperator{\cell}{{\brm{cell}}}
\DeclareMathOperator{\cof}{{\brm{cof}}}
\DeclareMathOperator{\inj}{{\brm{inj}}}

\DeclareMathOperator{\cyl}{{Cyl}}
\DeclareMathOperator{\CSA}{CSA}

\begin{document}

\begin{abstract} 
  We proved in a previous work that Cattani-Sassone's higher
  dimensional transition systems can be interpreted as a
  small-orthogonality class of a topological locally finitely
  presentable category of weak higher dimensional transition
  systems. In this paper, we turn our attention to the full
  subcategory of weak higher dimensional transition systems which are
  unions of cubes. It is proved that there exists a left proper
  combinatorial model structure such that two objects are weakly
  equivalent if and only if they have the same cubes after
  simplification of the labelling. This model structure is obtained by
  Bousfield localizing a model structure which is left determined with
  respect to a class of maps which is not the class of monomorphisms.
  We prove that the higher dimensional transition systems
  corresponding to two process algebras are weakly equivalent if and
  only if they are isomorphic. We also construct a second Bousfield
  localization in which two bisimilar cubical transition systems are
  weakly equivalent.  The appendix contains a technical lemma about
  smallness of weak factorization systems in coreflective
  subcategories which can be of independent interest.  This paper is a
  first step towards a homotopical interpretation of bisimulation for
  higher dimensional transition systems.
\end{abstract}

\maketitle

\tableofcontents

\section{Introduction}

\subsection*{Presentation of the paper}
Directed homotopy is a field of research aiming at studying the link
between concurrency and algebraic topology. In such a setting,
concurrency is modelled by higher-dimensional ``structures'' between
execution paths. In topological models like the ones of $d$-space
\cite{mg}, $d$-space generated by cubes \cite{FR}, flow \cite{model3},
globular complex \cite{diCW}, local po-space \cite{MR1683333}, locally
preordered space \cite{SK}, multipointed $d$-space
\cite{interpretation-glob}, these homotopies are homotopies in the
usual sense which preserve the direction of time. In combinatorial
models coming from the notion of (pre)cubical sets \cite{labelled}
\cite{exHDA} \cite{EWDCooperating} \cite{Pratt} \cite{Gunawardena1}
\cite{rvg} \cite{ccsprecub} \cite{symcub}, the concurrent execution of
$n$ actions is modelled by an $n$-cube, in which each axis of
coordinates corresponds to one action.

Concurrency is modelled in a somewhat different way in the formalism of
higher dimensional transition systems introduced by Cattani and
Sassone \cite{MR1461821}. Indeed, the concurrent execution of $n$
actions is modelled by a \emph{multiset} of $n$ actions. A multiset is
a set with possible repetition of some elements
(e.g. $\{0,0,2,3,3,3\}$). This notion is a generalization of the
$1$-dimensional notion of transition system in which transitions
between states are labelled by one action (e.g.,
\cite[Section~2.1]{MR1365754}). The latter $1$-dimensional notion
cannot of course model concurrency. It is proved in \cite{hdts} that
Cattani-Sassone's higher dimensional transition systems are a
small-orthogonality class of a larger category of \emph{weak higher
  dimensional transition systems (weak HDTS)} enjoying very nice
categorical properties: topological and locally finitely
presentable. Cattani-Sassone's higher dimensional transition systems
are weak HDTS satisfying two axioms CSA1 (cf. Definition~\ref{csa1})
and CSA2 (understood first and second Cattani-Sassone Axiom):
cf. Definition~\ref{isa} for a weaker form of CSA2. In plain English,
the first one says that one action between two given states can be
realized by at most one transition~\footnote{In CCS, the transition
  $a.P \stackrel{a} \rightarrow P$ is the unique transition from $a.P$
  to $P$.} The axiom CSA1 used by Cattani and Sassone is even stronger
(see the remark after Definition~\ref{csa1}) but we do not need it by
now. The second one is an analogue of the face operators in the
setting of precubical sets. These two axioms are satisfied by all
examples coming from process algebras.

It is not really a surprise that most of the topological models of
directed homotopy can be endowed with mathematical structures which
are very close to the ones existing in algebraic topology. In
particular, various model category structures can be related to
directed homotopy. It is more surprising that this kind of structure
exists in the setting of higher dimensional transition systems as
well.

We introduce in this paper the full subcategory of \emph{cubical
  transition systems}. A cubical transition system is a weak HDTS
which is equal to the union of its subcubes. Cubical transition
systems have a straightforward interpretation in concurrency. All
examples coming from process algebras are cubical because all these
examples are already colimits of cubes.  However, a cubical transition
system is not necessarily a colimit of cubes and the full subcategory
of weak HDTS generated by the colimits of cubes does not enjoy the
closure property we expect to find in such a setting. For example, the
boundary of the $2$-cube (cf. Definition~\ref{boundary-def}) is never
a colimit of cubes, but is always cubical.

The main result of this paper is that the category of cubical
transition systems can be endowed with a structure of left determined
left proper combinatorial model category structure with respect to a
class of cofibrations which is not the class of monomorphisms. This
model category structure is really minimal. Indeed, the corresponding
homotopy category cannot even identify all pairs of cubical transition
systems containing the same cubes ! We prove that there exists a
Bousfield localization such that two cubical transition systems are
weakly equivalent if and only if they have the same cubes after
simplification of the labelling. We also prove the existence of a
Bousfield localization with respect to the proper class of
bisimulations so that in the latter localization, two bisimilar
cubical transition systems are weakly equivalent.

\subsection*{Organization of the paper}

This paper starts in Section~\ref{remind-whdts} with a reminder about
weak higher dimensional transition systems (weak HDTS). Some
information about locally presentable and topological categories are
also collected here. It is important to say that the topological
structure plays an important role in the work, as well as the theory
of locally presentable categories which is extensively used, in
particular in Appendix~\ref{restriction}. Possible references for
these subjects are \cite{MR95j:18001} \cite{topologicalcat}
\cite{MR2506258} \cite{MR99h:55031}.

In Section~\ref{cubical-hdts}, we want to introduce the notion of
cubical transition system. Two equivalent definitions of them are
given: the weak HDTS equal to the union of their subcubes or
coreflective small-injectivity class. The last characterization
already implies that the category is locally presentable. It is
actually proved that it is locally \emph{finitely} presentable. It is
not topological since the adjunction between cubical transition
systems and weak HDTS is not concrete. Indeed, the coreflector removes
every action which is not used in a transition
(cf. Proposition~\ref{coreflector-explicit}). So what plays the role
of the underlying set varies. It is important to understand that the
full subcategory of cubes is not a dense or even a strong generator of
the category of cubical transition systems. It is necessary to add a
new family of weak HDTS, the \emph{double transition $\dd{x}$ labelled
  by $x$} for $x$ running over the set $\Sigma$ of labels
(cf. Definition~\ref{double-transition}).

Section~\ref{remind-comb} is a reminder about combinatorial model
categories, that is cofibrantly generated model categories
\cite{ref_model2} \cite{MR99h:55031} such that the underlying category
is locally presentable. Olschok's paper \cite{MO}, which generalizes
to locally presentable categories Cisinski's techniques for
constructing homotopical structures on toposes \cite{MR1924082}, plays
a fundamental role in this work. The notions of Grothendieck localizer
and of left determined model category are also recalled in this
section.

Section~\ref{cstr} expounds the construction of the combinatorial
model structure on weak HDTS. This model category carries a segment
object (which has nothing to do with the $1$-cube !) which is the key
to verifying all hypotheses of Olschok's theorems. This model category
is left proper since all objects are cofibrant. It is also left
determined with respect to its class of cofibrations, i.e. it is the
one with the smallest class of weak equivalences with our class of
cofibrations.  This class of weak equivalences is actually really
small, as we will see. A cofibration of weak HDTS is by definition a
map which is one-to-one on actions, but not necessarily on states. So
a map like $R:\{0,1\} \rightarrow \{0\}$ (a set being identified with
the weak HDTS with same set of states, no actions and no transitions)
is a cofibration of weak HDTS, and also of cubical transition systems
since every set is cubical as a disjoint sum of $0$-cubes. A similar
cofibration $R:\{0,1\} \rightarrow \{0\}$ exists in the model category
of flows \cite{model3} but we do not know whether there is a deeper
connexion between these two facts.

Section~\ref{homotopy-cube} restricts the previous structure to the
full subcategory of cubical transition systems. By definition, a
cofibration of cubical transition systems is a map between cubical
transition systems which is a cofibration of weak HDTS. The main
problem is to prove the smallness of the class of cofibrations between
cubical transition systems. The set of generating cofibrations used
for constructing the left determined model structure of $\whdts$
cannot be reused since they involve weak HDTS which are not
cubical. It is certainly possible to use combinatorial methods to find
a generating set of the class of cofibrations of cubical transition
systems. We use in this paper techniques of the theory of locally
presentable categories. This is the subject of
Appendix~\ref{restriction} which is of independent interest
(cf. Theorem~\ref{catlem}). The argument is a kind of generalization
of Smith's arguments to prove his well-known theorem
(Theorem~\ref{method1}), and more specifically for proving the
smallness of the class of trivial cofibrations. But let us repeat:
here the purpose is the proof of the smallness of the class of
cofibrations. The smallness of the class of trivial cofibrations is a
consequence of Olschok's theorems. This model category is also left
proper since all objects are cofibrant. It is also left determined
with respect to its class of cofibrations.

The next Section~\ref{weakisiso} characterizes the weak equivalences
in the left determined model structure of cubical transition
systems. It appears that CSA1 has a homotopical
interpretation. Roughly speaking, two cubical transition systems are
weakly equivalent in the left determined model structure if and only
if they are isomorphic modulo the first Cattani-Sassone axiom. It
follows that the canonical map $C_1[x] \sqcup C_1[x] \longrightarrow
\dd{x}$ sending two copies of the $1$-cube generated by $x$ to the
double transition labelled by $x$ is \emph{not} a weak equivalence
(cf. Figure~\ref{contrex-mono}). It is also proved in this section as
intermediate result that every cubical transition system which
satisfies CSA1 is fibrant.

Section~\ref{loc} overcomes this problem by proving that it is
possible to Bousfield localize with respect to the cubification
functor. The above map becomes a weak equivalence since $C_1[x] \sqcup
C_1[x]$ is precisely the cubification of $\dd{x}$.  In this Bousfield
localization, two cubical transition systems are weakly equivalent if
and only if they have the same cubes after simplification of the
labelling.

Finally Section~\ref{bisim} sketches the link with bisimulation. This
will be the subject of future works. 

Appendix~\ref{restriction} is the categorical lemma used in the core
of the paper which is of independent interest.

There are some remarks scattered in the paper about process algebras
with references to \cite{hdts}. But no knowledge about them is
required to read this paper and these remarks can be skipped without
problem.

\section{Weak higher dimensional transition systems}
\label{remind-whdts}

All categories are locally small. The set of maps in a category $\K$
from $X$ to $Y$ is denoted by $\K(X,Y)$. The locally small category
those objects are the maps of $\K$ and those morphisms are the
commutative squares is denoted by $\Mor(\K)$. The initial (final
resp.) object, if it exists, is always denoted by $\varnothing$
($\mathbf{1}$). The identity of an object $X$ is denoted by $\id_X$.
A subcategory will be by convention always \emph{isomorphism-closed}.

\begin{nota} A non empty set of {\rm labels} $\Sigma$ is fixed.  \end{nota}

Let us recall in this section the definition of a weak HDTS and some
fundamental examples. We start by collecting some well-known facts
about locally presentable and topological categories.

\subsection*{Locally presentable categories} Let $\lambda$ be a
regular cardinal, i.e. such that the poset $\lambda$ is
$\lambda$-directed \cite[p 160]{MR1697766}. An object $X$ of a
category $\K$ is \emph{$\lambda$-presentable} if the functor $\K(X,-)$
preserves $\lambda$-directed colimits.  A category $\K$ is
\emph{$\lambda$-accessible} if there exists a set of
$\lambda$-presentable objects such that every object of $\K$ is a
$\lambda$-directed colimit of objects of this set. A category $\K$ is
\emph{locally $\lambda$-presentable} if it is cocomplete and
$\lambda$-accessible. A subcategory $\LL$ of a category $\K$ is
\emph{accessibly-embedded} if it is full and closed under
$\lambda$-directed colimits for some regular cardinal $\lambda$. A
functor $F: \C \rightarrow \D$ is \emph{accessible} if there exists a
regular cardinal $\lambda$ such that $\C$ and $\D$ are
$\lambda$-accessible and $F$ preserves $\lambda$-directed colimits.
Every accessible functor satisfies the solution-set condition by
\cite[Corollary~2.45]{MR95j:18001}. When $\lambda = \aleph_0$, the
prefix ``$\lambda$-'' is replaced by ``finitely''.  In the preceding
definitions, $\lambda$-directed diagrams can be substituted by
$\lambda$-filtered diagrams by \cite[Remark~1.21]{MR95j:18001} since
for every (small) $\lambda$-filtered category $\D$, there exists a
(small) $\lambda$-directed poset $\D_0$ and a cofinal functor $\D_0
\rightarrow \D$.

\subsection*{Topological categories} The paradigm of \emph{topological
  category} over the category of $\set$ is the one of general
topological spaces with the notions of initial topology and final
topology \cite{topologicalcat}. More precisely, a functor $\omega: \C
\rightarrow \D$ is \emph{topological} (or $\C$ is \emph{topological}
over $\D$) if each cone $(f_i:X \rightarrow \omega A_i)_{i\in I}$
where $I$ is a class has a unique $\omega$-initial lift (the
\emph{initial structure}) $(\overline{f}_i:A \rightarrow A_i)_{i\in
  I}$, i.e.: 1) $\omega A= X$ and $\omega \overline{f}_i = f_i$ for
each $i\in I$; 2) given $h: \omega B \rightarrow X$ with $f_i h=
\omega \overline{h}_i$, $\overline{h}_i:B\rightarrow A_i$ for each
$i\in I$, then $h=\omega \overline{h}$ for a unique $\overline{h} : B
\rightarrow A$.  Topological functors can be characterized as functors
such that each cocone $(f_i:\omega A_i \rightarrow X)_{i\in I}$ where
$I$ is a class has a unique $\omega$-final lift (the \emph{final
  structure}) $\overline{f}_i:A_i \rightarrow A$, i.e.: 1) $\omega A=
X$ and $\omega \overline{f}_i = f_i$ for each $i\in I$; 2) given $h: X
\rightarrow \omega B$ with $hf_i = \omega \overline{h}_i$,
$\overline{h}_i:A_i \rightarrow B$ for each $i\in I$, then $h=\omega
\overline{h}$ for a unique $\overline{h} : A \rightarrow B$. Let us
suppose $\D$ complete and cocomplete. A limit (resp. colimit) in $\C$
is calculated by taking the limit (resp. colimit) in $\D$, and by
endowing it with the initial (resp. final) structure. In this work, a
topological category is a topological category over the category
$\set^{\{s\} \cup \Sigma}$ where $\{s\} \cup \Sigma$ is called the set
of sorts.

\subsection*{Weak higher dimensional transition systems (weak HDTS)}

\bd A {\rm weak higher dimensional transition system (weak HDTS)}
consists of a triple \[(S,\mu:L\rightarrow \Sigma,T=\bigcup_{n\geq
  1}T_n)\] where $S$ is a set of {\rm states}, where $L$ is a set of
{\rm actions}, where $\mu:L\rightarrow \Sigma$ is a set map called the
{\rm labelling map}, and finally where $T_n\subset S\p L^n\p S$ for $n
\geq 1$ is a set of {\rm $n$-transitions} or {\rm $n$-dimensional
  transitions} such that one has:
\begin{itemize}
\item (Multiset axiom) For every permutation $\sigma$ of
  $\{1,\dots,n\}$ with $n\geq 2$, if $(\alpha,u_1,\dots,u_n,\beta)$ is
  a transition, then $(\alpha,u_{\sigma(1)}, \dots, u_{\sigma(n)},
  \beta)$ is a transition as well.
\item (Coherence axiom) For every $(n+2)$-tuple
  $(\alpha,u_1,\dots,u_n,\beta)$ with $n\geq 3$, for every $p,q\geq 1$
  with $p+q<n$, if the five tuples $(\alpha,u_1, \dots, u_n, \beta)$,
  $(\alpha,u_1, \dots, u_p, \nu_1)$, $(\nu_1, u_{p+1}, \dots, u_n,
  \beta)$, $(\alpha, u_1, \dots, u_{p+q}, \nu_2)$ and $(\nu_2,
  u_{p+q+1}, \dots, u_n, \beta)$ are transitions, then the
  $(q+2)$-tuple $(\nu_1, u_{p+1}, \dots, u_{p+q}, \nu_2)$ is a
  transition as well.
\end{itemize}
A map of weak higher dimensional transition systems
\[f:(S,\mu : L \rightarrow \Sigma,(T_n)_{n\geq 1}) \rightarrow
(S',\mu' : L' \rightarrow \Sigma ,(T'_n)_{n\geq 1})\] consists of a
set map $f_0: S \rightarrow S'$, a commutative square
\[
\xymatrix{
  L \ar@{->}[r]^-{\mu} \ar@{->}[d]_-{\widetilde{f}}& \Sigma \ar@{=}[d]\\
  L' \ar@{->}[r]_-{\mu'} & \Sigma}
\] 
such that if $(\alpha,u_1,\dots,u_n,\beta)$ is a transition, then
$(f_0(\alpha),\widetilde{f}(u_1),\dots,\widetilde{f}(u_n),f_0(\beta))$
is a transition. The corresponding category is denoted by $\whdts$.
The $n$-transition $(\alpha,u_1,\dots,u_n,\beta)$ is also called a
{\rm transition from $\alpha$ to $\beta$}.  \ed

\begin{nota} The labelling map from the set of actions to the set of
  labels will be very often denoted by $\mu$. \end{nota}

A transition $(\alpha,u_1, \dots, u_n, \beta)$ intuitively means that
one goes from the state $\alpha$ to the state $\beta$ by executing
concurrently $n$ actions $u_1,\dots,u_n$. Hence the Multiset axiom,
which replaces the multiset formalism of \cite{MR1461821}. The
Coherence axiom is more complicated to understand. We just want to say
here that it is the topological part (in the sense of topological
categories) of an axiom introduced by Cattani and Sassone themselves
and that it is necessary for the mathematical development of the
theory: it is necessary to view Cattani-Sassone's higher dimensional
transition systems as a small-orthogonality class of $\whdts$. All
cubes satisfy this axiom and inside a given cube, the Coherence axiom
ensures that all transitions glue together properly.  Formally, this
axiom looks like a $5$-ary composition, even if it is topological.  We
refer to \cite{hdts} for further explanations.

The category $\whdts$ is locally finitely presentable by
\cite[Theorem~3.4]{hdts}. The functor \[\omega : \whdts
\longrightarrow \set^{\{s\}\cup \Sigma}\] taking the weak higher
dimensional transition system $(S,\mu : L \rightarrow
\Sigma,(T_n)_{n\geq 1})$ to the $(\{s\}\cup \Sigma)$-tuple of sets
$(S,(\mu^{-1}(x))_{x\in \Sigma}) \in \set^{\{s\}\cup \Sigma}$ is
topological by \cite[Theorem~3.4]{hdts} too.

\begin{nota} For $n\geq 1$, let $0_n = (0,\dots,0)$ ($n$-times) and
  $1_n = (1,\dots,1)$ ($n$-times). By convention, let 
  $0_0=1_0=()$. \end{nota}

We give now some important examples of weak HDTS. In each of the
following examples, the Multiset axiom and the Coherence axiom are
satisfied for trivial reasons. 

\begin{enumerate}
\item Let $n \geq 0$. Let $x_1,\dots,x_n \in \Sigma$. The \emph{pure
    $n$-transition} $C_n[x_1,\dots,x_n]^{ext}$ is the weak HDTS with
  the set of states $\{0_n,1_n\}$, with the set of actions $\{(x_1,1),
  \dots, (x_n,n)\}$ and with the transitions all $(n+2)$-tuples
  $(0_n,(x_{\sigma(1)},\sigma(1)), \dots,
  (x_{\sigma(n)},\sigma(n)),1_n)$ for $\sigma$ running over the set of
  permutations of the set $\{1,\dots ,n\}$.
\item Every set $X$ may be identified with the weak HDTS having the
  set of states $X$, with no actions and no transitions.
\item For every $x\in \Sigma$, let us denote by $\underline{x}$ the
  weak HDTS with no states, one action $x$, and no
  transitions. Warning: the weak HDTS $\{x\}$ contains one state $x$
  and no actions whereas the weak HDTS $\underline{x}$ contains no
  states and one action $x$.
\item For every $x\in \Sigma$, let us denote by $\dd{x}$ the weak HDTS
  with four states $\{1,2,3,4\}$, one action $x$ and two transitions
  $(1,x,2)$ and $(3,x,4)$.
\end{enumerate}

\bd \label{double-transition} The weak HDTS $\dd{x}$ is called the
{\rm double transition (labelled by $x$)} where $x\in \Sigma$.  \ed

Let us introduce now the weak HDTS corresponding to the $n$-cube.

\bp \label{cas_cube} \cite[Proposition~5.2]{hdts} Let $n\geq 0$ and
$x_1,\dots,x_n\in \Sigma$. Let $T_d\subset \{0,1\}^n \p
\{(x_1,1),\dots,(x_n,n)\}^d \p \{0,1\}^n$ (with $d\geq 1$) be the
subset of $(d+2)$-tuples
\[((\epsilon_1,\dots,\epsilon_n), (x_{i_1},i_1),\dots,(x_{i_d},i_d),
(\epsilon'_1,\dots,\epsilon'_n))\] such that
\begin{itemize}
\item $i_m = i_n$ implies $m = n$, i.e. there are no repetitions in the
  list $(x_{i_1},i_1),\dots,(x_{i_d},i_d)$
\item for all $i$, $\epsilon_i\leq \epsilon'_i$
\item $\epsilon_i\neq \epsilon'_i$ if and only if
  $i\in\{i_1,\dots,i_d\}$. 
\end{itemize}
Let $\mu : \{(x_1,1),\dots,(x_n,n)\} \rightarrow \Sigma$ be the set
map defined by $\mu(x_i,i) = x_i$. Then \[C_n[x_1,\dots,x_n] =
(\{0,1\}^n,\mu : \{(x_1,1),\dots,(x_n,n)\}\rightarrow
\Sigma,(T_d)_{d\geq 1})\] is a well-defined weak HDTS called the {\rm
  $n$-cube}. \ep

For $n = 0$, $C_0[]$, also denoted by $C_0$, is nothing else but the
weak HDTS $(\{()\},\mu:\varnothing \rightarrow
\Sigma,\varnothing)$. For every $x\in \Sigma$, one has $C_1[x] =
C_1[x]^{ext}$.  In \cite{hdts}, it is explained how the $n$-cube
$C_n[x_1,\dots,x_n]$ is freely generated by the pure $n$-transition
$C_n[x_1,\dots,x_n]^{ext}$. It is not necessary to recall this point
here.

\section{Cubical transition systems}
\label{cubical-hdts}

\subsection*{Definition of $\cts$} Before giving the definition of a
cubical transition system, we need first to check out that unions of
objects exist in $\whdts$.  So this section starts by studying the
monomorphisms of $\whdts$.

\bp \label{car_mono} A map $f:X = (S,\mu:L \rightarrow \Sigma,T)
\rightarrow X' = (S',\mu': L' \rightarrow \Sigma,T')$ of $\whdts$ is a
monomorphism if and only if the set maps $f_0: S \rightarrow S'$ and
$\widetilde{f}: L \rightarrow L'$ are one-to-one.  \ep

\bpf \underline{Only if part}. Suppose that $f : X \rightarrow X'$ is
a monomorphism. Let $\alpha$ and $\beta$ be two states of $X$ with
$f_0(\alpha) = f_0(\beta)$. Consider the two maps of weak higher
dimensional transition systems $g,h:\{0\} \rightarrow X$ defined by
$g(0) = \alpha$ and $h(0) = \beta$. Since $f$ is a monomorphism, one
has $g=h$. Therefore $\alpha=\beta$. Thus, the set map $f_0:S
\rightarrow S'$ is one-to-one. Now let $u$ and $v$ be two actions of
$X$ with $\widetilde{f}(u) = \widetilde{f}(v)$. One necessarily has
$\mu(u) = \mu(v) = x \in \Sigma$. Let $g,h : \underline{x} \rightarrow
X$ be the two maps of higher dimensional transition systems defined
respectively by $g(x) = u$ and $h(x) = v$. Then $g = h$ since $f$ is a
monomorphism. Therefore $u = v$ and $\widetilde{f}$ is one-to-one.
\underline{If part}. Let $f:X \rightarrow Y$ be a weak higher
dimensional transition system such that $f_0$ and $\widetilde{f}$ are
both one-to-one. Let $g,h:Z \rightarrow X$ be two maps of higher
dimensional transition systems such that $fg=fh$. Then $f_0g_0 =
f_0h_0$ and $\widetilde{f}\widetilde{g} =
\widetilde{f}\widetilde{h}$. So $g_0 = h_0$ and $\widetilde{g} =
\widetilde{h}$. The forgetful functor $\whdts \rightarrow \set^{\{s\}
  \cup \Sigma}$ is topological, and therefore faithful by
\cite[Theorem~21.3]{topologicalcat}. So $g = h$ and $f$ is a
monomorphism.  \epf

\bp Every family of subobjects of a weak HDTS has an union, i.e. a
least upper bound in the family of subobjects.  \ep

\bpf Let $(f_i : X_i \rightarrow X)_{i\in I}$ be a family of
subobjects of a weak HDTS $X$. Let $X_i = (S_i,\mu:L_i \rightarrow
\Sigma,T_i)$. Consider the set of states $S' = \bigcup_{i\in I}
(f_i)_0(S_i)$ and the set of actions $L' = \bigcup_{i\in I}
\widetilde{f_i}(L_i)$ equipped with the final structure. We obtain a
weak HDTS $X'$ and by Proposition~\ref{car_mono}, the canonical map
$X'\rightarrow X$ is a monomorphism. The weak HDTS $X'$ is the union
of the $(f_i : X_i \rightarrow X)_{i\in I}$.  \epf

We are now ready to give the definition of a cubical transition
system.

\bd Let $X$ be a weak HDTS. A {\rm cube} of $X$ is a map
$C_n[x_1,\dots,x_n] \longrightarrow X$. A {\rm subcube} of $X$ is the
image of a cube of $X$. A weak HDTS is a {\rm cubical transition
  system} if it is equal to the union of its subcubes. The full
subcategory of cubical transition systems is denoted by $\cts$.  \ed

Let $x_1,\dots,x_n \in \Sigma$ with $n\geq 0$. For $n\geq 2$, the weak
HDTS $C_n[x_1,\dots,x_n]^{ext}$ is not cubical since the union of its
subcubes is equal to its set of states $\{0_n,1_n\}$. The weak HDTS
$C_n[x_1,\dots,x_n]$ is always a cubical transition system since the
image of the identity of $C_n[x_1,\dots,x_n]$ is a subcube.  The weak
HDTS $\dd{x}$ is cubical for every $x\in \Sigma$. The weak HDTS
$\underline{x}$ is \emph{never} cubical for any $x\in \Sigma$ since
the union of its subcube is equal to $\varnothing$. For every set $A$,
the corresponding weak HDTS $A$ is cubical as a disjoint sum of
$0$-cubes.

\subsection*{Lifting property and small-injectivity class}

\bd Let $i:A\longrightarrow B$ and $p:X\longrightarrow Y$ be maps of
$\K$. Then $i$ has the {\rm left lifting property} (LLP) with respect
to $p$ (or $p$ has the {\rm right lifting property} (RLP) with respect
to $i$) if for every commutative square of solid arrows
\[
\xymatrix{
A\ar@{->}[dd]_{i} \ar@{->}[rr]^{\alpha} && X \ar@{->}[dd]^{p} \\
&&\\
B \ar@{-->}[rruu]^{k}\ar@{->}[rr]_{\beta} && Y,}
\]
there exists a morphism $k$ called a \textit{lift} making both
triangles commutative. This situation is denoted by $f \square g$. \ed

Let us introduce the notations $\inj_\K(\C) = \{g \in \K, \forall f
\in \C, f\square g\}$ and $\cof_\K(\C) = \{f \in \K, \forall g\in
\inj_\K(\C), f\square g\}$ where $\C$ is a class of maps of $\K$. The
class of morphisms of $\K$ that are transfinite compositions of
pushouts of elements of $\C$ is denoted by $\cell_\K(\C)$.  An element
of $\cell_\K(\C)$ is called a \textit{relative $\C$-cell complex}. The
cocompleteness of $\K$ implies $\cell_\K(\C)\subset \cof_\K(\C)$. When
the class $\C$ is a set $I$, every morphism of $\cof_\K(I)$ is a
retract of a morphism of $\cell_\K(I)$ by
\cite[Corollary~2.1.15]{MR99h:55031} since in a locally presentable
category, the domains of $I$ are always small relative to
$\cell_\K(I)$.

Sometimes, the letter $\K$ in the notations $\cof_\K$, $\inj_\K$ and
$\cell_\K$ may be omitted if the underlying category we are working
with is obvious.

\begin{center}
\underline{By convention, the letter $\K$ will be
  always omitted if $\K = \whdts$.}
\end{center}

\bd \cite[Definition~4.1]{MR95j:18001} Let $S$ be a set of
  maps of a locally presentable category $\K$. The full subcategory of
  $S$-injective objects (called a {\rm small-injectivity class}) of
  $\K$ is generated by $\{X \in \K \mid X\rightarrow \mathbf{1} \in
  \inj(S)\}$.
\ed

Let us recall that an object $X$ is \emph{orthogonal} to $S$ if not
only it is injective, but also the factorization is unique. A
small-injectivity class of a locally presentable category is always
accessible. A small-orthogonality class (the subclass of objects
orthogonal to a given set of objects) of a locally presentable
category is always a reflective locally presentable subcategory. Read
\cite[Chapter~1.C]{MR95j:18001} and \cite[Chapter~4]{MR95j:18001} for
further details. For an epimorphism $f$, being $f$-orthogonal is
equivalent to being $f$-injective.

\subsection*{The cubical transition systems as a small-injectivity class}

\bth \label{preacc} The category of cubical transition systems is a
small-injectivity class of $\whdts$.  More precisely, a weak HDTS $X$
is a cubical transition system if and only if it is injective with
respect to the set of inclusions $C_n[x_1,\dots,x_n]^{ext} \subset
C_n[x_1,\dots,x_n]$ and $\underline{x_1} \subset C_1[x_1]$ for all
$n\geq 0$ and all $x_1,\dots, x_n \in \Sigma$. \eth

\bpf \underline{Only if part}. 1) Let $X$ be a cubical transition
system.  Let $C_n[x_1,\dots,x_n]^{ext} \rightarrow X$ be a map of weak
HDTS.  Let $(\alpha,u_1,\dots,u_n,\beta)$ be the image by this map of
the transition $(0_n,(x_1,1),\dots,(x_n,n),1_n)$.  By hypothesis,
there exists a cube $C_m[y_1,\dots,y_m] \rightarrow X$ of $X$ such
that the image contains the transition
$(\alpha,u_1,\dots,u_n,\beta)$. There is not yet any reason for $m$ to
be equal to $n$. This means that the image of $C_m[y_1,\dots,y_m]
\rightarrow X$ contains the image of $C_n[x_1,\dots,x_n]^{ext}
\rightarrow X$. In other terms, the latter map factors as a composite
\[C_n[x_1,\dots,x_n]^{ext} \longrightarrow C_m[y_1,\dots,y_m]
\longrightarrow X.\] By \cite[Theorem~5.6]{hdts}, the map
$C_n[x_1,\dots,x_n]^{ext} \rightarrow C_m[y_1,\dots,y_m]$ factors as a
composite $C_n[x_1,\dots,x_n]^{ext} \rightarrow C_n[x_1,\dots,x_n]
\rightarrow C_m[y_1,\dots,y_m]$ since the cube $C_m[y_1,\dots,y_m]$ is
injective, and even orthogonal to the inclusion
$C_n[x_1,\dots,x_n]^{ext} \subset
C_n[x_1,\dots,x_n]$~\footnote{Orthogonality means that this
  factorization is unique but we do not need this fact here.}.  Thus,
$X$ is injective with respect to the set of maps
$C_n[x_1,\dots,x_n]^{ext} \subset C_n[x_1,\dots,x_n]$ for all $n\geq
0$ and all $x_1,\dots, x_n \in \Sigma$.  2) Let $\underline{x_1}
\rightarrow X$ be a map of weak HDTS. By hypothesis, there exists a
cube $C_m[y_1,\dots,y_m] \rightarrow X$ of $X$ such that the image
contains the image of $\underline{x_1} \rightarrow X$. In other terms,
the latter map factors as a composite
\[\underline{x_1} \longrightarrow C_m[y_1,\dots,y_m] \longrightarrow
X.\] Since the maps of weak HDTS preserve labellings, there exists $k$
such that $x_1 = y_k$.  Hence the factorization
\[\underline{x_1} \longrightarrow C_1[x_1] \longrightarrow
C_m[y_1,\dots,y_m] \longrightarrow X.\] So $X$ is injective with
respect to the set of maps $\underline{x_1} \subset C_1[x_1]$ for
$x_1$ running over $\Sigma$. \underline{If part}. Every transition and
every state of $X$ belong to a subcube since $X$ is injective with
respect to the maps $C_n[x_1,\dots,x_n]^{ext} \subset
C_n[x_1,\dots,x_n]$ for all $n\geq 0$ and all $x_1,\dots, x_n \in
\Sigma$. Every action of $X$ belongs to a subcube because $X$ is
injective with respect to the maps $\underline{x_1} \subset C_1[x_1]$
for $x_1$ running over $\Sigma$.  \epf

It follows that the category $\cts$ of cubical transition systems is
accessible by \cite[Proposition~4.7]{MR95j:18001}. It is even locally
finitely presentable, as we will see.

\subsection*{Some elementary facts about (co)reflective subcategories}

A \emph{coreflective (resp. reflective) subcategory} of a category
$\C$ is a full isomorphism-closed category such that the inclusion
functor is a left (resp. right) adjoint. The right (resp. left)
adjoint is called the \emph{coreflector} (resp. the \emph{reflector}).
The two following propositions are elementary and well-known. We use
them several times so we need to state them clearly.

\bp \label{facile_1} \cite[page~89]{MR1712872} Let $\D\subset \C$ be a
coreflective (isomorphism-closed) subcategory of a category $\C$,
i.e. a full subcategory such that the inclusion $\D \subset \C$ has a
right adjoint $R:\C \rightarrow \D$. Then:
\begin{enumerate}
\item The counit $R(X) \rightarrow X$ is an isomorphism if and only if
  $X$ belongs to $\D$
\item If $\C$ is cocomplete, then so is $\D$. 
\end{enumerate} \ep

\bp \label{colim_gen} \cite[Proposition~3.1(i)]{1185.18014} Let $\C$
be a cocomplete category. Let $\mathcal{S}$ be a set of objects of
$\C$. The full subcategory of colimits of objects of $\mathcal{S}$ is
a coreflective subcategory $\C_\mathcal{S}$ of $\C$.  The right
adjoint to the inclusion functor $\C_\mathcal{S} \subset \C$ is the
``Kelleyfication'' functor $k_\mathcal{S}$ defined by:
\[k_\mathcal{S}(X) = \liminj_{\begin{array}{c}S\rightarrow X\\S\in
    \mathcal{S}\end{array}} S.\] 
\ep

\subsection*{Coreflectivity of the category of cubical transition systems}

First we recall how colimits are calculated in $\whdts$.

\bp \label{calcul_colim} \cite[Proposition~3.5]{hdts} Let $X = \liminj
X_i$ be a colimit of weak higher dimensional transition systems with
$X_i = (S_i,\mu_i:L_i \rightarrow \Sigma, T^i = \bigcup_{n\geq 1}
T^i_n)$ and $X = (S,\mu:L \rightarrow \Sigma, T = \bigcup_{n\geq 1}
T_n)$. Then:
\begin{enumerate}
\item $S = \liminj S_i$, $L = \liminj L_i$, $\mu = \liminj \mu_i$
\item the union $\bigcup_{i} T^i$ of the image of the $T^i$ in
  $\bigcup_{n\geq 1} (S \p L^n\p S)$ satisfies the Multiset axiom.
\item $T$ is the closure of $\bigcup_{i} T^i$ under the Coherence
  axiom. 
\item when the union $\bigcup_{i} T^i$ is already closed under the
  Coherence axiom, this union is the final structure.
\end{enumerate}
\ep 

\begin{lem} \label{used}
  Consider a colimit $\liminj X_i$ in $\whdts$ such that every action
  $u$ of $X_i$ is used, i.e. there exists a transition
  $(\alpha_i,u_i,\beta_i)$ of $X_i$. Then every action of $X$ is used.
\end{lem}

\bpf By Proposition~\ref{calcul_colim}, the set of transitions of
$\liminj X_i$ is obtained by taking the closure under the Coherence
axiom of the union of the transitions of the $X_i$, hence the result
since the set of actions of $\liminj X_i$ is the union of the actions
of the $X_i$. \epf

\bth \label{right-adjoint}
Let $X\in \whdts$. The counit map 
\[q_X : \liminj_{{\begin{array}{c}f: C_n[x_1,\dots,x_n]
    \rightarrow X\\\hbox{ or }f:\dd{x} \rightarrow X\end{array}}}
\dom(f) \rightarrow X\] where $\dom(f)$ is the domain of $f$ is
bijective on states and one-to-one on actions and transitions.
Moreover, the weak HDTS $X$ is cubical if and only if $q_X$ is an
isomorphism.  \eth

\bpf It is important to keep in mind that, since $\whdts$ is
topological, the set of states (resp. of actions) of $\dom(q_X)$ is
the colimit of the sets of states (resp. of actions) of the $\dom(f)$
for $f$ running over the set of maps of the form $C_n[x_1,\dots,x_n]
\rightarrow X$ or $\dd{x} \rightarrow X$ for $n\geq 0$,
$x_1,\dots,x_n,x\in \Sigma$. 

\underline{$q_X$ is one-to-one on states}. Let $\alpha$ and
$\beta$ be two states of $\dom(q_X)$ having the same image $\gamma$ in
$X$. Then the diagram $\{\alpha\} \leftarrow \{\gamma\} \rightarrow
\{\beta\}$ is a subdiagram in the colimit calculating $\dom(q_X)$. Hence
$\alpha = \gamma = \beta$ in $\dom(q_X)$.

\underline{$q_X$ is onto on states}. Let $\alpha$ be a state of
$X$. Then the map $C_0[] \rightarrow X$ mapping the unique state of
$C_0[]$ to $\alpha$ is in the colimit calculating
$\dom(q_X)$. 

\underline{$q_X$ is one-to-one on actions}. Let $u$ and $v$ be two
actions of $\dom(q_X)$ having the same image $w$ in $X$. By
Lemma~\ref{used}, the maps $\underline{u} \rightarrow \dom(q_X)$ and
$\underline{v} \rightarrow \dom(q_X)$ factor as
composites \[\underline{u} \longrightarrow C_1[\mu(u)] \longrightarrow
\dom(q_X) \hbox{ and }\underline{v} \longrightarrow C_1[\mu(v)]
\longrightarrow \dom(q_X).\] One has $\mu(u) = \mu(v) = \mu(w) = x \in
\Sigma$ by definition of a map of weak HDTS.  Therefore, there exists
a commutative diagram of weak HDTS like in Figure~\ref{crucial} Hence
$u=v$ in $\dom(q_X)$.
\begin{figure}
\[
\xymatrix{
 C_1[\mu(u)] \ar@{->}[rd] \ar@/^10pt/[rrrd]&& & \\
& \dd{x} \fR{}&&  X\\
C_1[\mu(v)]\ar@{->}[ru] \ar@/_10pt/[rrru] && &
}
\]
\caption{The crucial role of $\dd{x}$}
\label{crucial}
\end{figure}

\underline{$q_X$ is one-to-one on transitions}.  Let
$(\alpha,u_1,\dots,u_n,\beta)$ and
$(\alpha',u'_1,\dots,u'_{n'},\beta')$ be two transitions of
$\dom(q_X)$ having the same image in $X$. Then one has $n = n'$. Since
$q_X$ is one-to-one on states, one gets $\alpha = \alpha'$ and $\beta
= \beta'$. Since $q_X$ is one-to-one on actions, one gets $u_i = u'_i$
for $1\leq i \leq n$.

Let us prove now the last part of the theorem. Let $X$ be a cubical
transition system.  Let $u$ be an action of $X$. Then there exists a
map $\underline{\mu(u)} \rightarrow X$ mapping $\mu(u)$ to $u$. By
Theorem~\ref{preacc}, the latter map factors as a composite
\[\underline{\mu(u)} \longrightarrow C_1[\mu(u)] \longrightarrow
X\] since $X$ is cubical. Hence $q_X$ is onto on actions. Let
$(\alpha,u_1,\dots,u_n,\beta)$ be a transition of $X$. Then there
exists a map $C_n[\mu(u_1),\dots,\mu(u_n)]^{ext} \rightarrow X$
mapping the transition \[(0_n,(\mu(u_1),1),\dots,(\mu(u_n),n),1_n)\] to
$(\alpha,u_1,\dots,u_n,\beta)$.  By Theorem~\ref{preacc}, the latter
map factors as a composite \[C_n[\mu(u_1),\dots,\mu(u_n)]^{ext}
\longrightarrow C_n[\mu(u_1),\dots,\mu(u_n)] \longrightarrow X\] since
$X$ is cubical.  Hence $q_X$ is onto on transitions.  So $q_X$
is an isomorphism.  Conversely, let us suppose now that $q_X$ is an
isomorphism. Let $f: \underline{x} \rightarrow X$ be a map of weak
HDTS. Then, by hypothesis, the action $\widetilde{f}(x)$ of $X$ comes
from an action $u$ of $\dom(q_X)$. The corresponding map
$\underline{x} = \underline{\mu(u)} \rightarrow \dom(q_X)$ factors as
a composite
\[\underline{x} = \underline{\mu(u)} \longrightarrow C_1[\mu(u)]
\longrightarrow \dom(q_X)\] by construction of $q_X$. Hence $X$ is
injective with respect to the maps $\underline{x} \rightarrow C_1[x]$
for $x\in \Sigma$. Let $g : C_n[x_1,\dots,x_n]^{ext} \rightarrow X$ be
a map of weak HDTS. Then, by hypothesis, the transition
$(g_0(0_n),\widetilde{g}(x_1,1),\dots, \widetilde{g}(x_n,n),g_0(1_n))$
of $X$ comes from a transition $(\alpha,u_1,\dots,u_n,\beta)$ of
$\dom(q_X)$. The corresponding map $C_n[\mu(u_1),\dots,\mu(u_n)]^{ext}
\rightarrow \dom(q_X)$ factors as a composite
\[C_n[\mu(u_1),\dots,\mu(u_n)]^{ext} \longrightarrow
C_n[\mu(u_1),\dots,\mu(u_n)] \longrightarrow \dom(q_X)\] by
construction of $q_X$. Hence $X$ is injective with respect to the maps
\[C_n[\mu(u_1),\dots,\mu(u_n)]^{ext} \longrightarrow
C_n[\mu(u_1),\dots,\mu(u_n)].\] So by Theorem~\ref{preacc}, the weak
HDTS $X$ is cubical.  \epf

\begin{cor} The full subcategory of $\cts$ generated by the cubes
$C_n[x_1,\dots,x_n]$ for $n\geq 0$ and $x_1,\dots,x_n\in \Sigma$ and
by the weak HDTS $\dd{x}$ for $x\in \Sigma$ is dense in
$\cts$.  \end{cor}

\bd Let $X\in \whdts$. The {\rm cubification functor} is the functor
\[\cub : \whdts \longrightarrow \whdts\] defined by \[\cub =
\liminj_{C_n[x_1,\dots,x_n] \rightarrow X} C_n[x_1,\dots,x_n].\]
Denote by $p_X:\cub(X) \rightarrow X$ the canonical map.  \ed

The full subcategory generated by the cubes $C_n[x_1,\dots,x_n]$ for
$n\geq 0$ and $x_1,\dots,x_n\in \Sigma$ is not a dense, and even not a
strong generator of $\cts$. It is not a dense generator since the weak
HDTS $\dd{x}$ is not a colimit of cubes.  Indeed, the canonical
map \[C_1[x] \sqcup C_1[x] \iso \cub(\dd{x}) \longrightarrow \dd{x}\]
is not an isomorphism. The left-hand weak HDTS contains two distinct
actions $x_1$ and $x_2$ labelled by $x$, whereas the right-hand one
contains only one action $x$. It is not a strong generator either
since the canonical map (cf. Figure~\ref{contrex-mono})
\[\cub(\dd{x}) \longrightarrow \dd{x}\] is a monomorphism in $\cts$~\footnote{It is not a monomorphism in $\whdts$: the precompositions by $\underline{x} \rightarrow C_1[x] \sqcup C_1[x]$ mapping $x$ to $x_1$ and to $x_2$ give the same result.}
and since every map $C_n[x_1,\dots,x_n] \rightarrow \dd{x}$ factors as
a composite $C_n[x_1,\dots,x_n] \rightarrow C_1[x] \sqcup C_1[x]
\rightarrow \dd{x}$ ($n$ is necessarily equal to $1$).

\begin{rem} The map of Figure~\ref{contrex-mono} is also an
  epimorphism. 
\end{rem}

\begin{figure}
\[\left\{\begin{array}{c} C_1[x] \sqcup C_1[x]) \\
        {\stackrel{x_1}\longrightarrow} \\
        {\stackrel{x_2}\longrightarrow}\end{array}\right.
  \stackrel{p_{x}}\longrightarrow \left\{\begin{array}{c} \liminj (C_1[x] \leftarrow \underline{x}
        \rightarrow C_1[x]) \\ {\stackrel{x}\longrightarrow}\\
        {\stackrel{x}\longrightarrow}\end{array}\right.\]
\caption{Monomorphism in $\cts$ with $\mu(x_1) = \mu(x_2) = x$}
\label{contrex-mono}
\end{figure}

\begin{cor} The category $\cts$ is a coreflective locally finitely
presentable subcategory of $\whdts$. \end{cor}

\bpf The right adjoint to the inclusion functor $\cts \subset \whdts$
is the functor $X\mapsto \dom(q_X)$ by
Proposition~\ref{colim_gen}. The category is therefore cocomplete with
set of dense (and therefore strong) finitely presentable generators
the cubes $C_n[x_1,\dots,x_n]$ for $n\geq 0$ and $x_1,\dots,x_n\in
\Sigma$ and the weak HDTS $\dd{x}$ for $x\in \Sigma$. The category
$\cts$ is therefore locally finitely presentable by
\cite[Theorem~1.20]{MR95j:18001}.  \epf

\bd \label{boundary-def} Let $n\geq 1$ and $x_1,\dots,x_n \in
\Sigma$. Let $\de C_n[x_1,\dots,x_n]$ be the weak HDTS defined by
removing from its set of transitions all $n$-transitions. It is called
the {\rm boundary} of $C_n[x_1,\dots,x_n]$. \ed

The weak HDTS $\de C_2[x_1,x_2]$ is not a colimit of cubes but is
cubical: it is obtained by identifying states in the cubical
transition system $\dd{x_1} \sqcup \dd{x_2}$.

\section{About combinatorial model categories}
\label{remind-comb}

\bd \cite{MR2003h:18001} Let $\K$ be a locally presentable category. A
{\rm weak factorization system} is a pair $(\mathcal{L},\mathcal{R})$
of classes of morphisms of $\K$ such that $\inj_\K(\mathcal{L}) =
\mathcal{R}$ and such that every morphism of $\K$ factors as a
composite $r\circ \ell$ with $\ell\in \mathcal{L}$ and $r\in
\mathcal{R}$. The weak factorization system is {\rm functorial} if the
factorization $r\circ \ell$ can be made functorial.  \ed

For every set of maps $I$ of a locally presentable category $\K$, the
pair of classes of maps $(\cof_\K(I), \inj_\K(I))$ is a weak
factorization system by \cite[Proposition~1.3]{MR1780498}. A weak
factorization system of the form $(\cof_\K(I), \inj_\K(I))$ is said
\textit{small}, or {\rm generated by $I$}.  A small weak factorization
system is necessarily functorial.

For every weak factorization system $(\mathcal{L},\mathcal{R})$, the
class of maps $\mathcal{L}$ is closed under retract, pushout and
transfinite composition.

\bd \cite{MR99h:55031} A {\rm combinatorial model category} is a
locally presentable category equipped with three classes of morphisms
$(\C,\F,\W)$ (resp. called the classes of {\rm cofibrations}, {\rm
fibrations} and {\rm weak equivalences}) such that:
\begin{enumerate}
\item the class of morphisms $\W$ is closed under retracts and
satisfies the two-out-of-three axiom i.e.: if $f$ and $g$ are
morphisms of $\K$ such that $g\circ f$ is defined and two of $f$, $g$
and $g\circ f$ are weak equivalences, then so is the third.
\item the pairs $(\C\cap \W, \F)$ and $(\C, \F\cap \W)$ are both small
weak factorization systems. So there exist two sets of maps $I$ and
$J$ such that $(\C, \F\cap \W) = (\cof_\K(I), \inj_\K(I))$ and
$(\C\cap \W, \F) = (\cof_\K(J), \inj_\K(J))$.
\end{enumerate} The triple $(\C,\F,\W)$ is called a {\rm model
category structure}. An element of $\C \cap \W$ is called a {\rm
trivial cofibration}. An element of $\F\cap \W$ is called a {\rm
trivial fibration}. A map of $I$ is called a {\rm generating
cofibration} and a map of $J$ a {\rm generating trivial
cofibration}. \ed

There exists at most one model category structure $(\C,\F,\W)$ for a
given class of cofibrations $\C$ and a given class of weak
equivalences $\W$. Indeed, the class of cofibrations determines the
class of trivial fibrations, and the intersection of the classes of
cofibrations and of weak equivalences determines the class of
fibrations.

An object $X$ is \emph{cofibrant} (\emph{fibrant} resp.) if the
canonical map $\varnothing \rightarrow X$ ($X \rightarrow \mathbf{1}$)
is a cofibration (fibration resp.).  A model category is
\emph{left proper} if the pushout along a cofibration of a weak
equivalence is a weak equivalence. By a well-know theorem due to
C. L. Reedy \cite{Reedy}, every model category such that every object
is cofibrant is left proper (e.g.,
\cite[Corollary~13.1.3]{ref_model2}).

For every object $X$ of a model category, the canonical map
$\varnothing \rightarrow X$ ($X \rightarrow \mathbf{1}$ resp.) factors
as a composite $0 \rightarrow X^{cof} \rightarrow X$ ($X \rightarrow
X^{fib} \rightarrow \mathbf{1}$ resp.) where $X^{cof}$ is cofibrant
and $X^{cof} \rightarrow X$ is a trivial fibration ($X^{fib}$ is
fibrant and $X \rightarrow X^{fib}$ is a trivial cofibration
resp.). $X^{cof}$ ($X^{fib}$ resp.) is called the \emph{cofibrant (fibrant
resp.) replacement functor}.

\bd \cite[Definition~3.4]{MR1924082} Let $\mathcal{A}$ be a class of
morphisms of a category $\K$. A class of maps $\W$ satisfying the
two-out-of-three axiom, such that $\inj_\K(\mathcal{A}) \subset \W$
and such that $\mathcal{A} \cap \W$ is closed under pushout and
transfinite composition is called a {\rm $\mathcal{A}$-localizer}, or
a {\rm localizer with respect to $\mathcal{A}$}. \ed

The class of all maps of $\K$ is clearly an $\mathcal{A}$-localizer
and the intersection of any family of $\mathcal{A}$-localizers is a
$\mathcal{A}$-localizer. Therefore there exists a smallest
$\mathcal{A}$-localizer containing a given set of maps $S$ denoted by
$\W^\K_\mathcal{A}(S)$, or $\W_\mathcal{A}(S)$ if there is no
ambiguity (once again, $\K$ will be always omitted if $\K = \whdts$).

Let $\K$ be a locally presentable category. Let $\mathcal{A}$ be a
class of morphisms of $\K$. There exists at most one model structure
on $\K$ such that $\mathcal{A}$ is the class of cofibrations and such
that $\W_\mathcal{A}(\varnothing)$ is the class of weak equivalences
since the class of trivial cofibrations is then completely known and
by definition of a weak factorization system, the classes of
fibrations and trivial fibrations are determined as well.  When it
exists, it is called the \emph{left determined model structure with
  respect to $\mathcal{A}$} \cite{rotho}. Note that the existence of
this model structure implies that $\W_\mathcal{A}(\varnothing)$ is
closed under retract. However, this hypothesis is not in the
definition of a localizer.

\bd \cite{ideeloc} A {\rm very good cylinder} of a weak factorization
system $(\mathcal{L},\mathcal{R})$ in a locally presentable category
$\K$ is a functorial factorization of the codiagonal $X\sqcup X
\rightarrow X$ as a composite \[\xymatrix@1{ X\sqcup X \fR{\gamma_X}
  && {\cyl(X)} \fR{\sigma_X}&& X}\] with $\gamma_X \in \mathcal{L}$
and $\sigma_X \in \mathcal{R}$.  Two maps $f,g:X \rightrightarrows Y$
are {\rm homotopy equivalent} if the pair $(f,g)$ belongs to the
symmetric transitive closure of the binary relation $f\sim g$ whenever
the map $f\sqcup g:X\sqcup X \rightarrow Y$ factors as a composite
\[\xymatrix@1{ X\sqcup X \fR{\gamma_X} && {\cyl(X)} \fR{H}&& Y.}\] 
The homotopy relation does not depend on the choice of a very good
cylinder by \cite[Observation~3.3]{ideeloc}.  \ed

The adjective \emph{very good} (meaning that $\sigma_X \in
\mathcal{R}$) is not used in \cite{ideeloc}. The adjective
\emph{final} is used in \cite{MO}. The terminology of
\cite[Definition~4.2]{MR1361887} seems to be better to avoid any
confusion with the notion of final structure in a topological
category.

\begin{nota} The two composites 
\[\xymatrix@1{ X \subset X\sqcup X \fR{\gamma_X} && {\cyl(X)}}\]
are denoted by $\gamma_X^0$ and $\gamma_X^1$.
\end{nota}

\begin{nota} For every map $f:X \rightarrow Y$ and every natural
  transformation $\alpha : F \Rightarrow F'$ between two endofunctors
  of $\K$, the map $f\star \alpha$ is the canonical map
\[f\star \alpha : FY \sqcup_{FX} F'X \longrightarrow F'Y\] 
induced by the commutative diagram of solid arrows 
\[\xymatrix{
FX \fR{\alpha_X}\fD{Ff} && F'X \fD{F'f} \\
&& \\
FY \fR{\alpha_Y} && F'Y}
\]
and the universal property of the pushout.
\end{nota}

\bd \cite[Definition~3.8]{MO} A very good cylinder of a weak
factorization system $(\mathcal{L},\mathcal{R})$ in a locally
presentable category $\K$ is {\rm cartesian} if the cylinder functor
$\cyl : \K \rightarrow \K$ is a left adjoint and if one has the
inclusions $\mathcal{L} \star \gamma \subset \mathcal{L}$ and
$\mathcal{L} \star \gamma^k \subset \mathcal{L}$ for $k=0,1$. \ed

A \emph{cylinder of a model category} is a very good cylinder for the
weak factorization system formed by the cofibrations and the trivial
fibrations.

Let us conclude the section by recalling well-known Smith's theorem
generating model structures on locally presentable categories.

\bth (Smith) \label{method1} Let $I$ be a set of morphisms of a
locally presentable category $\K$. Let $\W$ be an accessible
accessibly-embedded $\cof_\K(I)$-localizer closed under retracts. Then
there exists a cofibrantly generated model structure on $\K$ with
class of cofibrations $\cof_\K(I)$, with class of fibrations
$\inj_\K(\cof_\K(I) \cap \W)$, and with class of weak equivalences
$\W$.  \eth

\bpf[Sketch of proof] The class $\W$ satisfies the solution set
condition by \cite[Corollary~2.45]{MR95j:18001}. Hence the existence
of the model structure by Smith's theorem
\cite[Theorem~1.7]{MR1780498}. \epf

The \emph{Bousfield localization} of a model category $\mathcal{M}$ by
a class of maps $\LL$ is a model category $L_\LL \mathcal{M}$ with the
same underlying category, the same class of cofibrations, together with a
map of model categories~\footnote{i.e. a left adjoint preserving
  cofibrations and trivial cofibrations} $\mathcal{M} \rightarrow
L_\LL \mathcal{M}$ such that every map of model categories
$\mathcal{M} \rightarrow \mathcal{N}$ taking the cofibrant replacement
of every map of $\LL$ to a weak equivalence of $\mathcal{N}$ factors
uniquely as a composite $\mathcal{M} \rightarrow L_\LL \mathcal{M}
\rightarrow \mathcal{N}$. The properties of this object used in this
paper are listed now:
\begin{enumerate}
\item The Bousfield localization of a left proper combinatorial model
  category with respect to any \emph{set} of maps always exists and is
  left proper combinatorial \cite{MR2506258} \cite{Lurie}
  \cite[Theorem~3.3.19]{ref_model2}.
\item A weak equivalence between two cofibrant-fibrant objects in
  $L_\LL \mathcal{M}$ is a weak equivalence of $\mathcal{M}$
  \cite[Theorem~3.2.13]{ref_model2}.
\end{enumerate}
By Bousfield localization of $\mathcal{M}$ with respect to a functor
$F : \mathcal{M} \rightarrow \mathcal{M}$ preserving weak
equivalences, it is meant the Bousfield localization with respect to the
class of maps $f$ such that $F(f)$ is a weak equivalence.

\section{The left determined model category of weak HDTS}
\label{cstr}

The purpose of this section is the proof of the existence of the
left determined model structure with respect to the cofibrations of
weak HDTS defined as follows:

\bd A {\rm cofibration} of weak HDTS is a map of weak HDTS inducing an
injection between the set of actions.  \ed

Note that the class of cofibrations is strictly bigger than the class
of monomorphisms of $\whdts$ since $R:\{0,1\} \rightarrow \{0\}$ is a
cofibration of weak HDTS. We do not know if there is a link between
this fact and the existence of an analogous cofibration on the model
category of flows introduced in \cite{model3}.

\bp \label{lem_prelim} The class of cofibrations of weak HDTS is
closed under pushout, transfinite composition and retract. \ep

\bpf Since the functor $\omega : \whdts \longrightarrow
\set^{\{s\}\cup \Sigma}$ is topological, it is colimit-preserving. So
it suffices to observe that the class of injections in the category of
sets is closed under retract, pushout and transfinite composition, for
example by considering the weak factorization system of the category
of sets $(\cof_\set(C),\inj_\set(C))$ where $C:\varnothing \subset
\{0\}$ denotes the inclusion.  \epf

\begin{nota} Let $\I$ be the set of maps $C:\varnothing \rightarrow
  \{0\}$, $R:\{0,1\} \rightarrow \{0\}$, $\varnothing \subset
  \underline{x}$ for $x\in \Sigma$ and $\{0_n,1_n\} \sqcup
  \underline{x_1} \sqcup \dots \sqcup \underline{x_n} \subset C_n[x_1,
  \dots, x_n]^{ext}$ for $n\geq 1$ and $x_1, \dots, x_n \in
  \Sigma$. \end{nota}

\bp\label{small} One has $\cell(\I) = \cof(\I)$ and this class of maps
is the class of cofibrations of weak HDTS.  \ep

\bpf Every map of $\I$ is a cofibration of weak HDTS. Since $\I$ is a
set, the class of maps $\cof(\I)$ is the closure under retract of
transfinite composition of pushouts of elements of $\I$. So $\cell(\I)
\subset \cof(\I)$ and by Proposition~\ref{lem_prelim}, every map of
$\cof(\I)$ is a cofibration of weak HDTS. It then suffices to prove
that every cofibration of weak HDTS belongs to $\cell(\I)$.

Let $f : X = (S, \mu : L \rightarrow \Sigma, T) \rightarrow X' = (S',
\mu': L' \rightarrow \Sigma, T')$ be a cofibration of weak HDTS. The
set map $f_0:S \rightarrow S'$ factors as a composite $S \rightarrow
f_0(S) \subset S'$. The left-hand map is a transfinite composition of
pushouts of $R:\{0,1\} \rightarrow \{0\}$. The inclusion $f_0(S)
\subset S'$ is a transfinite composition of pushouts of $C:\varnothing
\rightarrow \{0\}$. By hypothesis, the set map $\widetilde{f} : L
\rightarrow L'$ is one-to-one. Consider the pushout diagram of weak
HDTS
\[
\xymatrix{S \sqcup \left(\bigsqcup_{u\in L}\underline{\mu(u)}\right)  \fR{\subset} \fD{f \sqcup \widetilde{f}} && X \fD{} \\
&&\\
S' \sqcup \left(\bigsqcup_{u\in L'}\underline{\mu'(u)}\right) \fR{} && \cocartesien Y. }
\] 
The universal property of the pushout yields 
a map of weak HDTS $g : Y \rightarrow X'$ such that $g_0$ and
$\widetilde{g}$ are bijections.  Consider the pushout diagram of weak
HDTS
\[
\xymatrix{\bigsqcup\limits_{(\alpha,u_1,\dots,u_n,\beta)\in T'\backslash T} (\{0_n,1_n\}\sqcup
\underline{\mu'(u_1)} \sqcup \dots \sqcup \underline{\mu'(u_n)})
\fRr{{\begin{array}{c}0_n\mapsto \alpha\\1_n\mapsto \beta\\ \mu'(u_i)
    \mapsto \mu'(u_i)\end{array}}}\fD{}
&&& Y \fD{}\\
&&& \\
\bigsqcup\limits_{(\alpha,u_1,\dots,u_n,\beta)\in T'\backslash T} C_n[\mu'(u_1), \dots,
\mu'(u_n)]^{ext} \fRr{}   &&&  \cocartesien Z.}\]
The universal property of the pushout yields 
a map
$ h: Z \rightarrow X'$ such that $h_0$ and $\widetilde{h}$ are
bijections. So the set of transitions of $Z$ can be identified with a
subset of the set of transitions of $X'$. By construction, the map $h$
induces an onto map between the set of transitions. So $h$ is an
isomorphism of weak HDTS and $\cell(\I) = \cof(\I)$. \epf

The terminal object $\mathbf{1}$ of $\whdts$ is described as follows:
the set of states is $\{0\}$, the set of actions is $\Sigma$, the
labelling map is the identity of $\Sigma$ and the set of transitions
is $\bigcup_{n \geq 1} \Sigma^n$. In other terms, one has $\mathbf{1}
\iso (\{0\}, \id_\Sigma, \bigcup_{n \geq 1} \Sigma^n)$. Let $V$ be the
weak HDTS
\[\boxed{V := (\{0\}, \pr_1 : \Sigma\p \{0,1\} \rightarrow \Sigma,
  \{0\} \p (\bigcup_{n \geq 1} (\Sigma\p \{0,1\})^n) \p \{0\})}\]
$V$ is called the \textit{segment object} of $\whdts$. 

\bp \label{produit_hdts} Let $X = (S,\mu:L \rightarrow \Sigma, T)$ and
$X' = (S',\mu':L' \rightarrow \Sigma, T')$ be two weak HDTS. The
binary product $X\p X'$ has the set of states $S\p S'$, the set of
actions $L\p_\Sigma L' = \{(x,x')\in L\p L', \mu(x) = \mu'(x')\}$ and
the labelling map $\mu\p_\Sigma \mu':L\p_\Sigma L' \rightarrow
\Sigma$. A tuple $((\alpha,\alpha'), (u_1,u'_1), \dots, (u_n,u'_n),
(\beta,\beta'))$ is a transition of $X\p X'$ if and only if $\mu(u_i)
= \mu'(u'_i)$ for $1\leq i \leq n$ with $n\geq 1$, the tuple $(\alpha,
u_1, \dots, u_n, \beta)$ is a transition of $X$ and $(\alpha', u'_1,
\dots, u'_n, \beta')$ a transition of $X'$. \ep

\bpf The forgetful functor $\omega : \whdts \longrightarrow
\set^{\{s\}\cup \Sigma}$ is limit-preserving by
\cite[Proposition~21.12]{topologicalcat} since it is topological. So
the set of states is $S\p S'$, the set of actions $L\p_\Sigma L'$ and
the labelling map $\mu\p_\Sigma \mu':L\p_\Sigma L' \rightarrow
\Sigma$. Consider the set $T'''$ of tuples $((\alpha,\alpha'),
(u_1,u'_1), \dots, (u_n,u'_n), (\beta,\beta'))$ such that $\mu(u_i) =
\mu'(u'_i)$ for $1\leq i \leq n$ with $n\geq 1$, the tuple $(\alpha,
u_1, \dots, u_n, \beta)$ is a transition of $X$ and $(\alpha', u'_1,
\dots, u'_n, \beta')$ a transition of $X'$. The existence of the
projections $X\p X' \rightarrow X$ and $X\p X' \rightarrow X'$ implies
that the set of transitions $T''$ of $X \p X'$ satisfies $T'' \subset
T'''$.  Let $t = (\alpha, u_1, \dots, u_n, \beta) \in T$ and $t' =
(\alpha', u'_1, \dots, u'_n, \beta') \in T'$ such that $\mu(u_i) =
\mu'(u'_i)$ for $1\leq i \leq n$ with $n\geq 1$.  Let $t\p t'$ be the
weak HDTS with set of states $S\p S'$, with set of actions $L\p_\Sigma
L'$, with labelling map $\mu\p_\Sigma \mu'$, and with set of
transitions $\{((\alpha,\alpha'), (u_{\sigma(1)},u'_{\sigma(1)}),
\dots, (u_{\sigma(n)},u'_{\sigma(n)}), (\beta,\beta')), \sigma \hbox{
  permutation of }\{1,\dots,n\}\}$. Since the set of transitions $T''$
is given by an initial structure, the cone of weak HDTS $(t\p t'
\rightarrow X, t\p t' \rightarrow X')$ induced by the projections
factors uniquely by a map $t\p t' \rightarrow X\p X'$ which is the
identity on the set of states and the set of actions.  So $T'''\subset
T''$.  \epf

\bp \label{somme_hdts} Let $X = (S,\mu:L \rightarrow \Sigma, T)$ and
$X' = (S',\mu':L' \rightarrow \Sigma, T')$ be two weak higher
dimensional transition systems. The binary coproduct $X\sqcup X'$ has
the set of states $S\sqcup S'$, the set of actions $L\sqcup L'$ and
the labelling map $\mu\sqcup \mu':L\sqcup L' \rightarrow \Sigma$. A
tuple $(\alpha, u_1, \dots, u_n, \beta)$ is a transition of $X\sqcup
X'$ if and only if it is a transition of $X$ or a transition of
$X'$. \ep

\bpf The forgetful functor $\omega : \whdts \longrightarrow
\set^{\{s\}\cup \Sigma}$ is colimit-preserving by
\cite[Proposition~21.12]{topologicalcat} since it is topological. So
the set of states is $S\sqcup S'$, the set of actions $L\sqcup L'$ and
the labelling map $\mu\sqcup \mu':L\sqcup L' \rightarrow \Sigma$.  The
disjoint union of the transitions of $X$ and $X'$ is closed under the
Coherence axiom. So it is equal to the set of transitions of $X\sqcup
X'$ by Proposition~\ref{calcul_colim}.  \epf

\bp \label{rlp} The canonical map $\mathbf{1} \sqcup \mathbf{1}
\rightarrow \mathbf{1}$ factors as a composite 
$\mathbf{1} \sqcup \mathbf{1} \longrightarrow V \longrightarrow
\mathbf{1}$ such that the left-hand map is a cofibration and such that
the right-hand map satisfies the right lifting property with respect
to every cofibration. \ep

\bpf Proposition~\ref{somme_hdts} tells us that the set of states
(resp. of actions) of $\mathbf{1} \sqcup \mathbf{1}$ is the disjoint
union of the set of states (resp. of actions) of $\mathbf{1}$.  Let
$\mathbf{1} \sqcup \mathbf{1} \rightarrow V$ be the map of weak HDTS
defined on states by the constant set map ($V$ has only one state) and
on actions by the bijection $\Sigma \sqcup \Sigma \rightarrow \Sigma
\p \{0,1\}$ taking the left-hand copy $\Sigma$ to $\Sigma \p \{0\}$
and the right-hand copy of $\Sigma$ to $\Sigma \p \{1\}$.  The
composite $\mathbf{1} \sqcup \mathbf{1} \rightarrow V\rightarrow
\mathbf{1}$ is the unique map of weak HDTS from $\mathbf{1} \sqcup
\mathbf{1}$ to $\mathbf{1}$. The map $\mathbf{1} \sqcup \mathbf{1}
\rightarrow V$ is a cofibration.

Consider the commutative square of solid arrows
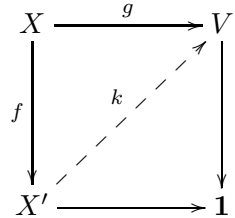
\[
\xymatrix{
X \fR{g} \fD{f} && V \fD{}\\
&& \\
X' \fR{} \ar@{-->}[rruu]^-{k} && \mathbf{1}}
\] 
where $f : X \rightarrow X'$ is a cofibration of weak HDTS. Let $X =
(S,\mu : L \rightarrow \Sigma, T)$ and $X' = (S',\mu' : L' \rightarrow
\Sigma, T')$.  Since $V$ has only one state, the definition of $k_0$
is clear: $k_0=0$. Since $f$ is a cofibration, $L$ can be identified
with a subset of $L'$.  Let $\widetilde{k} : L' \rightarrow \Sigma\p
\{0,1\}$ be the set map defined as follows:
\begin{itemize}
\item $\widetilde{k}(u) = \widetilde{g}(u)$ if $u \in L$ (we have no choice here)
\item $\widetilde{k}(u) = (\mu'(u),0)$ if $u \in L'\backslash L$.
\end{itemize}
Let $(\alpha,u_1,\dots,u_n,\beta)$ be a transition of $X'$. One always
has $\widetilde{k}(u_i) \in \{\mu'(u_i)\} \p \{0,1\}$, and necessarily
$\widetilde{k}(u_i) = (\mu'(u_i),0)$ if $u_i \in L'\backslash L$ for
every $i\in \{1,\dots,n\}$.  So the set maps $k_0$ and $\widetilde{k}$
takes the transition $(\alpha,u_1,\dots,u_n,\beta)$ to the tuple $(0,
(\mu'(u_1),\epsilon_1), \dots, (\mu'(u_n),\epsilon_n), 0)$ with
$\epsilon_1, \dots, \epsilon_{n} \in \{0,1\}$. The tuple $(0,
(\mu'(u_1),\epsilon_1), \dots, (\mu'(u_n),\epsilon_n), 0)$ is a
transition of $V$ by definition of $V$. So $k$ is a map of weak HDTS
and the map $V \rightarrow \mathbf{1}$ satisfies the RLP with respect
to every cofibration.  \epf

\bp \label{exp} The weak HDTS $V$ is exponentiable, i.e. the functor
$V\p - : \whdts \rightarrow \whdts$ has a right adjoint denoted by
$(-)^V: \whdts \rightarrow \whdts$. \ep

\bpf Let $Y = (S_Y, \mu:L_Y \rightarrow \Sigma, T_Y)$ be a weak HDTS. 
Recall that 
\[V := (\{0\}, \pr_1 : \Sigma\p \{0,1\} \rightarrow \Sigma, \{0\} \p
(\bigcup_{n \geq 1} (\Sigma\p \{0,1\})^n) \p \{0\}).\] Let us describe
at first the right adjoint
\[Y^V = (S^V,\mu^V: L^V \rightarrow \Sigma, T^V).\] One must have the
bijection of sets \[\whdts(V\p\{0\},Y) \iso \whdts(\{0\},Y^V) \iso
S^V.\] By Proposition~\ref{produit_hdts}, one has $V \p \{0\} \iso
\{0\}$. So necessarily there is the equality $S^V = S_Y$.  Let $x\in
\Sigma$. One must have the bijection of sets \[\whdts(V \p
\underline{x}, Y) \iso \whdts(\underline{x},Y^V) = (\mu^V)^{-1}(x).\]
By Proposition~\ref{produit_hdts} again, one has $V\p \underline{x}
\iso \underline{x} \sqcup \underline{x}$. Therefore one has
\[
(\mu^V)^{-1}(x) \iso \whdts(\underline{x}\sqcup\underline{x},X) \iso
\mu^{-1}(x)\p \mu^{-1}(x).\] Thus, one must necessarily have $L^V =
L_Y \p_\Sigma L_Y$ (the fibered product of $L_Y$ by itself over
$\Sigma$).  Finally, one must have the bijection of sets \[\whdts(V\p
C_n[x_1,\dots,x_n]^{ext},Y) \iso
\whdts(C_n[x_1,\dots,x_n]^{ext},Y^V)\] for every $x_1, \dots, x_n\in
\Sigma$.  By Proposition~\ref{produit_hdts} again, the $n$-transitions
of $Y^V$ are of the form $(\alpha,(u^-_1,u^+_1), \dots, (u^-_n,u^+_n),
\beta)$ such that the $2^n$ tuples $(\alpha, u_1^\pm, \dots, u_n^\pm,
\beta)$ are transitions of $Y$.

Let $X = (S_X, \mu:L_X \rightarrow \Sigma, T_X)$ be another weak HDTS.
Using Proposition~\ref{produit_hdts} again, let us describe now the
binary product $X\p V$. The set of states of $X\p V$ is $S_X$, the set
of actions is $L_X\p_\Sigma(\Sigma \p \{0,1\}) = L_X\p \{0,1\}$ and a
tuple $(\alpha,(u_1,\epsilon_1), \dots,(u_n,\epsilon_n),\beta)$ is a
transition if and only if $(\alpha,u_1, \dots,u_n,\beta)$ is a
transition of $X$.

The bijection $\whdts(X \p V,Y) \iso \whdts(X,Y^V)$ is then easy to
check.
\epf

\begin{nota} Let $\cyl(X) := X\p V$.
\end{nota}

\bp \label{fin_cartesien} One has $\cof(\I) \star \gamma^0 \subset
\cof(\I)$, $\cof(\I) \star \gamma^1 \subset \cof(\I)$ and $\cof(\I) \star
\gamma \subset \cof(\I)$. \ep

\bpf Let $f:X \rightarrow X'$ be a cofibration of weak HDTS. Let $X =
(S,\mu : L \rightarrow \Sigma, T)$ and $X' = (S',\mu' : L' \rightarrow
\Sigma, T')$.  The map of weak HDTS $f\star \gamma: (X'\sqcup
X')\sqcup_{X\sqcup X} \cyl(X) \rightarrow \cyl(X')$ is a cofibration
since the set map $\widetilde{f\star \gamma}$ is the identity of $L'
\sqcup L'$. The map of weak HDTS $f\star \gamma^k: X'\sqcup_{X}
\cyl(X) \rightarrow \cyl(X')$, where $\gamma^k_X: X\rightarrow
\cyl(X)$ and $\gamma^k_{X'}: X'\rightarrow \cyl(X')$ are the canonical
maps is a cofibration of weak HDTS since the set map $\widetilde{f
  \star \gamma^k}$ is the inclusion $L\sqcup L' \rightarrow L'\sqcup
L'$.  \epf

\bth \label{constr_model} Let $S$ be an arbitrary set of maps of
$\whdts$. The triple
\[(\cof(\I),\inj(\cof(\I)\cap \W_{\cof(\I)}(S)),\W_{\cof(\I)}(S))\] is
a left proper combinatorial model structure of $\whdts$. The segment
object $V$ is fibrant and contractible (i.e. weakly equivalent to the
terminal object) for this model structure. All objects are
cofibrant. \eth

\bpf By Proposition~\ref{exp}, Proposition~\ref{fin_cartesien} and
Proposition~\ref{rlp}, the functor $\cyl(X) = V\p X$ is a cartesian
very good cylinder for the weak factorization system
$(\cof(\I),\inj(\I))$. The latter weak factorization system is
cofibrant, i.e. all maps $\varnothing \rightarrow X$ belongs to
$\cof(\I)$ by Proposition~\ref{small}. The theorem is therefore a
consequence of \cite[Corollary~4.6]{MO}. \epf

When $S = \varnothing$, the above model structure is left determined
in the sense of \cite{rotho}, i.e. the class of weak equivalences is
the smallest localizer closed under retract. Indeed,
$\W_{\cof(\I)}(S)$ is included in this smallest localizer closed under
retract and it is closed under retract itself since it is the class of
weak equivalences of a model category structure.

Note that the category $\whdts$ is distributive in the following sense: 

\bp The category $\whdts$ is distributive, i.e. for every weak higher
dimensional transition system $X$, $Y$ and $Z$, there is the
isomorphism $(X\p Y) \sqcup (X \p Z) \iso X\p (Y \sqcup Z)$. \ep

\bpf Since the forgetful functor $\whdts \rightarrow \set^{\{s\}\cup
  \Sigma}$ is topological, it preserves limits and colimits by
\cite[Proposition~21.12]{topologicalcat}. So the canonical map $(X\p
Y) \sqcup (X \p Z) \rightarrow X\p (Y \sqcup Z)$ induces a bijection
between the sets of states and the sets of actions. So the set of
transitions $T$ of $(X\p Y) \sqcup (X \p Z)$ can be identified with a
subset of the set of transitions $T'$ of $X\p (Y \sqcup Z)$. So
$T\subset T'$. By Proposition~\ref{produit_hdts}, a transition of $X\p
(Y \sqcup Z)$ is of the form $((\alpha,\gamma), (u_1,v_1), \dots,
(u_n,v_n), (\beta,\delta))$ where the tuple $(\alpha,u_1, \dots, u_n,
\beta)$ is a transition of $X$ and where the tuple $(\gamma, v_1,
\dots, v_n, \delta)$ is a transition of $Y\sqcup Z$. By
Proposition~\ref{somme_hdts}, the transition $(\gamma, v_1, \dots,
v_n, \delta)$ is then either a transition of $Y$ or a transition of
$Z$. So by Proposition~\ref{produit_hdts} again, the tuple
$((\alpha,\gamma), (u_1,v_1), \dots, (u_n,v_n), (\beta,\delta))$ is
either a transition of $X\p Y$ or a transition of $X\p Z$. Thus,
$T'\subset T$.  \epf

The class of cofibrations is also stable under pullback along any map
(not necessarily product projection).  Therefore,
\cite[Remark~4.7]{MO} applies here: any factorization of the
codiagonal $\mathbf{1} + \mathbf{1} \rightarrow \mathbf{1}$ as a
composite $\mathbf{1} + \mathbf{1} \rightarrow W'
\rightarrow\mathbf{1}$ with the left-hand map a cofibration and the
right-hand map an element of $\inj(\I)$ will provide a very good
cylinder.

\section{The left determined model category of cubical transition systems}
\label{homotopy-cube}

In this section, $\mathcal{A}$ is a coreflective full subcategory of
$\whdts$.

\bth \label{constr_submodel} Let $\mathcal{A}$ be a coreflective
accessible subcategory of $\whdts$ such that:
\begin{itemize}
\item The class of cofibrations of $\whdts$ between objects of
  $\mathcal{A}$ is generated by a set, i.e. there exists a set
  $I_\mathcal{A}$ of maps of $\mathcal{A}$ such that
  $\cof_\mathcal{A}(I_\mathcal{A})$ is this class of maps.
\item The segment object $V$ belongs to $\mathcal{A}$.
\item The inclusion functor $\mathcal{A} \subset \whdts$ preserves
  binary products by $V$.
\end{itemize}
Let $S$ be an arbitrary set of maps of $\mathcal{A}$.
The triple
\[(\cof_\mathcal{A}(I_\mathcal{A}),\inj_\mathcal{A}(\cof_\mathcal{A}(I)\cap
\W^{\mathcal{A}}_{\cof(I_\mathcal{A})}(S)),\W^{\mathcal{A}}_{\cof(I_\mathcal{A})}(S))\]
is a left proper combinatorial model structure of $\mathcal{A}$.  \eth

\bpf The category $\mathcal{A}$ is cocomplete by
Proposition~\ref{facile_1}. Therefore it is locally presentable. So
the cylinder functor $X \mapsto V\p X$ is a left adjoint. The proof
then goes as for that of Theorem~\ref{constr_model}. The latter
theorem is in fact the particular case $\mathcal{A} = \whdts$. \epf

When $S = \varnothing$, the above model structure is left determined
in the sense of \cite{rotho}, i.e. the class of weak equivalences is
the smallest localizer closed under retract. 

\begin{nota} Let $\Lambda_\mathcal{A}(\cyl,S,I_\mathcal{A})$ be the set of maps:
\begin{itemize}
\item $\Lambda_\mathcal{A}^0(\cyl,S,I_\mathcal{A}) = S \cup (I_\mathcal{A} \star
  \gamma^0) \cup (I_\mathcal{A} \star \gamma^1)$
\item $\Lambda_\mathcal{A}^{n+1}(\cyl,S,I_\mathcal{A}) =
  \Lambda_\mathcal{A}^n(\cyl,S,I_\mathcal{A}) \star \gamma$
\item $\Lambda_\mathcal{A}(\cyl,S,I_\mathcal{A}) = \bigcup_{n\geq 0}
  \Lambda_\mathcal{A}^n(\cyl,S,I_\mathcal{A})$.
\end{itemize} 
\end{nota}

By \cite[Theorem~3.16, Theorem~4.5 and corollary~4.6]{MO}, the class
of weak equivalences $\W^{\mathcal{A}}_{\cof(I_\mathcal{A})}(S)$
coincides with the class of maps denoted by
$\W(\Lambda_\mathcal{A}(\cyl,S,I_\mathcal{A}))$ defined as follows. A
map $f : X \rightarrow Y$ of $\mathcal{A}$ belongs to
$\W(\Lambda_\mathcal{A}(\cyl,S,I_\mathcal{A}))$ if and only if for
every object $T$ of $\mathcal{A}$ such that the canonical map
$T\rightarrow \mathbf{1} \in
\inj_\mathcal{A}(\Lambda_\mathcal{A}(\cyl,S,I_\mathcal{A}))$, the
induced set map
\[\whdts(Y,T)/\simeq \longrightarrow \whdts(X,T)/\simeq \] 
is a bijection where $\simeq$ means the homotopy relation associated
with the cylinder $\cyl$. Moreover, the fibrant objects of the model
category of Theorem~\ref{constr_submodel} are exactly the objects $T$
such that $T\rightarrow \mathbf{1} \in
\inj_\mathcal{A}(\Lambda_\mathcal{A}(\cyl,S,I_\mathcal{A}))$.

\bth \label{constr_Quillen_adjoint} Let $\mathcal{A}$ and
$\mathcal{B}$ be two coreflective accessible subcategories of $\whdts$
with $\mathcal{A} \subset \mathcal{B}$ satisfying the hypotheses of
Theorem~\ref{constr_submodel}.  Let us suppose that the class of
cofibrations of $\whdts$ between objects of $\mathcal{A}$
(resp. $\mathcal{B}$) is generated by a set $I_\mathcal{A}$
(resp. $I_\mathcal{B}$). Let $S$ be an arbitrary set
of maps of $\mathcal{A}$. Let us equip $\mathcal{A}$ with the model structure 
\[(\cof_\mathcal{A}(I_\mathcal{A}),\inj_\mathcal{A}(\cof_\mathcal{A}(I)\cap
\W^{\mathcal{A}}_{\cof(I_\mathcal{A})}(S)),\W^{\mathcal{A}}_{\cof(I_\mathcal{A})}(S))\] and 
$\mathcal{B}$ with the model structure 
\[(\cof_\mathcal{B}(I_\mathcal{B}),\inj_\mathcal{B}(\cof_\mathcal{B}(I)\cap
\W^{\mathcal{B}}_{\cof(I_\mathcal{B})}(S)),\W^{\mathcal{B}}_{\cof(I_\mathcal{B})}(S)).\]
Then the inclusion functor $\mathcal{A} \subset \mathcal{B}$ is a left
Quillen adjoint.
\eth

\bpf The two categories $\mathcal{A}$ and $\mathcal{B}$ are cocomplete
by Proposition~\ref{facile_1} and therefore locally presentable. Since
the inclusion functor $\mathcal{A} \subset \mathcal{B}$ preserves
colimits (which are the same as the colimits of $\whdts$), it is a
left adjoint. it is clear that the inclusion functor takes
cofibrations to cofibrations. We must prove that it takes trivial
cofibrations to trivial cofibrations. It actually takes every weak
equivalence to a weak equivalence. Let $X\rightarrow Y$ be a weak
equivalence of $\mathcal{A}$.  Let $T$ be a fibrant object of
$\mathcal{B}$. Then the map $T \rightarrow \mathbf{1}$ satisfies the
RLP with respect to any map of
$\Lambda_\mathcal{A}(\cyl,S,I_\mathcal{A}) \subset
\Lambda_\mathcal{B}(\cyl,S,I_\mathcal{B})$.  So by adjunction, $R(T)
\rightarrow \mathbf{1}$ satisfies the RLP with respect to the maps of
$\Lambda_\mathcal{A}(\cyl,S,I_\mathcal{A})$, where $R(-)$ is the right
adjoint to the inclusion functor. So $R(T)$ is fibrant in
$\mathcal{A}$. Therefore the induced set map
\[\whdts(Y,R(T))/\simeq \longrightarrow \whdts(X,R(T))/\simeq \]
is a bijection.  So by adjunction again, $X\rightarrow Y$ is a weak
equivalence of $\mathcal{B}$. 
\epf

We want to apply Theorem~\ref{constr_submodel} to the case
$\mathcal{A} = \cts$ and $\mathcal{B} = \whdts$. 

\bd \label{isa} A weak HDTS $X$ satisfies the {\rm Intermediate state
  axiom} if for every $n\geq 2$, every $p$ with $1\leq p<n$ and every
transition $(\alpha,u_1,\dots,u_n,\beta)$ of $X$, there exists a (not
necessarily unique) state $\nu$ such that both
$(\alpha,u_1,\dots,u_p,\nu)$ and $(\nu,u_{p+1},\dots,u_n,\beta)$ are
transitions.  \ed

Note that the Unique intermediate state axiom CSA2 introduced in
\cite{hdts} is slightly stronger than the axiom above. Indeed, it
states that the intermediate states in a higher dimensional transition
are unique.

\bp \label{construction_map} \cite[Proposition~5.5]{hdts} Let $n\geq
0$ and $a_1,\dots,a_n\in \Sigma$. Let $X = (S,\mu:L \rightarrow
\Sigma, T = \bigcup_{n\geq 1} T_n)$ be a weak higher dimensional
transition system. Let $f_0: \{0,1\}^n \rightarrow S$ and
$\widetilde{f} : \{(a_1,1),\dots,(a_n,n)\} \rightarrow L$
be two set maps. Then the following conditions are equivalent:
\begin{enumerate}
\item The pair $(f_0,\widetilde{f})$ induces a map of weak higher
  dimensional transition systems from $C_n[a_1,\dots,a_n]$ to $X$.
\item For every transition $((\epsilon_1, \dots, \epsilon_n),
  (a_{i_1},i_1), \dots, (a_{i_r},i_r), (\epsilon'_1, \dots,
  \epsilon'_n))$ of $C_n[a_1,\dots,a_n]$ with $(\epsilon_1, \dots,
  \epsilon_n) = 0_n$ or $(\epsilon'_1, \dots, \epsilon'_n) = 1_n$, the
  tuple $(f_0(\epsilon_1, \dots, \epsilon_n),
  \widetilde{f}(a_{i_1},i_1), \dots, \widetilde{f}(a_{i_r}\linebreak[4],i_r),
  f_0(\epsilon'_1, \dots, \epsilon'_n))$ is a transition of $X$.
\end{enumerate} \ep

\bp \label{pre_charac_cub} A weak HDTS satisfies the Intermediate state
axiom if and only if it is injective with respect to the maps
$C_n[x_1,\dots,x_n]^{ext} \subset C_n[x_1,\dots,x_n]$ for all $n\geq
0$ and all $x_1,\dots, x_n \in \Sigma$. \ep

Recall that if a weak HDTS satisfies the Unique intermediate state
axiom CSA2, not only it is injective with respect to the maps
$C_n[x_1,\dots,x_n]^{ext} \subset C_n[x_1,\dots,x_n]$ for all $n\geq
0$ and all $x_1,\dots, x_n \in \Sigma$, but also the factorization is
unique: i.e. the weak HDTS is orthogonal to this set of maps
\cite[Theorem~5.6]{hdts}.

\bpf The proof is essentially an adaptation of the one of
\cite[Theorem~5.6]{hdts}. 

\underline{Only if part}. Let $X = (S,\mu:L \rightarrow \Sigma, T =
\bigcup_{n\geq 1} T_n)$ be a weak HDTS satisfying the Intermediate
state axiom.  Let $n\geq 0$ and $x_1,\dots,x_n \in \Sigma$. We have to
prove that the inclusion of weak HDTS $C_n[x_1,\dots,x_n]^{ext}
\subset C_n[x_1,\dots,x_n]$ induces an onto set map
\[\whdts(C_n[x_1,\dots,x_n],X) \longrightarrow
\whdts(C_n[x_1,\dots,x_n]^{ext},X).\] This fact is trivial for $n=0$
and $n=1$ since the inclusion $C_n[x_1,\dots,x_n]^{ext} \subset
C_n[x_1,\dots,x_n]$ is an equality. Let $f : C_n[x_1,\dots,x_n]^{ext}
\rightarrow X$ be a map of weak HDTS. The map $f$ induces a set map
$f_0:\{0_n,1_n\} \rightarrow S$ and a set map $\widetilde{f} :
\{(x_1,1),\dots,(x_n,n)\} \rightarrow L$.  Let
$(\epsilon_1,\dots,\epsilon_n) \in [n]$ be a state of
$C_n[x_1,\dots,x_n]$ different from $0_n$ and $1_n$.  Then there exist
(at least) two transitions \[(0_n,(x_{i_1},i_1),\dots,(x_{i_r},i_r),
(\epsilon_1,\dots,\epsilon_n))\] and
\[((\epsilon_1,\dots,\epsilon_n),(x_{i_{r+1}},i_{r+1}),\dots,(x_{i_{r+s}},i_{r+s}),
1_n)\] of $C_n[x_1,\dots,x_n]$ with $r,s\geq 1$. Let
$f_0(\epsilon_1,\dots,\epsilon_n)$ be a state of $X$ such that
\[(f_0(0_n),\widetilde{f}(x_{i_1},i_1),\dots,\widetilde{f}(x_{i_r},i_r),
f_0(\epsilon_1,\dots,\epsilon_n))\] and
\[(f_0(\epsilon_1,\dots,\epsilon_n),\widetilde{f}(x_{i_{r+1}},i_{r+1}),
\dots,\widetilde{f}(x_{i_{r+s}},i_{r+s}), f_0(1_n))\] are two
transitions of $X$. Since every transition from $0_n$ to
$(\epsilon_1,\dots,\epsilon_n)$ is of the form \[(0_n,
(x_{i_{\sigma(1)}},i_{\sigma(1)}), \dots,
(x_{i_{\sigma(r)}},i_{\sigma(r)}), (\epsilon_1,\dots,\epsilon_n))\]
where $\sigma$ is a permutation of $\{1,\dots,r\}$ and since every
transition from $(\epsilon_1,\dots,\epsilon_n)$ to $1_n$ is of the
form \[((\epsilon_1,\dots,\epsilon_n),
(x_{i_{\sigma'(r+1)}},i_{\sigma'(r+1)}), \dots,
(x_{i_{\sigma'(r+s)}},i_{\sigma'(r+s)}), 1_n)\] where $\sigma'$ is a
permutation of $\{r+1,\dots,r+s\}$, one obtains a well-defined set map
$f_0 : [n] \rightarrow S$.  The pair of set maps $(f_0,\widetilde{f})$
induces a well-defined map of weak HDTS by
Proposition~\ref{construction_map}. Therefore the set map
\[\whdts(C_n[x_1,\dots,x_n],X)  \longrightarrow
\whdts(C_n[x_1,\dots,x_n]^{ext},X)\] is onto.

\underline{If part}. Conversely, let $X = (S,\mu:L \rightarrow \Sigma,
T = \bigcup_{n\geq 1} T_n)$ be a weak HDTS injective to the set of
inclusions $\{C_n[x_1,\dots,x_n]^{ext} \subset C_n[x_1,\dots,x_n],
n\geq 0 \hbox{ and }x_1,\dots,x_n\in\Sigma\}$. Let
$(\alpha,u_1,\dots,u_n,\beta)$ be a transition of $X$ with $n\geq
2$. Then there exists a (unique) map
$C_n[\mu(u_1),\dots,\mu(u_n)]^{ext} \longrightarrow X$ taking the
transition $(0_n,(\mu(u_1),1),\dots,(\mu(u_n),n),1_n)$ to the
transition $(\alpha,u_1,\dots,u_n,\beta)$. By hypothesis, this map
factors as a composite \[C_n[\mu(u_1),\dots,\mu(u_n)]^{ext} \subset
C_n[\mu(u_1),\dots,\mu(u_n)] \stackrel{g}\longrightarrow X.\] Let
$1\leq p<n$. There exists a (unique) state $\nu$ of
$C_n[\mu(u_1),\dots,\mu(u_n)]$ such that the tuples $(0_n,
(\mu(u_1),1), \dots, (\mu(u_p),p), \nu)$ and $(\nu,
(\mu(u_{p+1}),p+1), \dots, (\mu(u_n),n), 1_n)$ are two transitions of
the HDTS $C_n[\mu(u_1), \dots, \mu(u_n)]$ by
Proposition~\ref{cas_cube}. Hence the existence of a state $g_0(\nu)$
of $X$ such that the tuples $(\alpha, u_1, \dots, u_p, g_0(\nu))$ and
$(g_0(\nu), u_{p+1}, \dots, u_n, \beta)$ are two transitions of
$X$. Thus, the weak HDTS $X$ satisfies the Intermediate state axiom.
\epf

\bp \label{charac_cub} A weak HDTS is a cubical transition system if
and only if it satisfies the Intermediate state axiom and every action
$u$ is {\rm used} in at least one $1$-transition $(\alpha,u,\beta)$. 
\ep

\bpf The statement is a corollary of Proposition~\ref{pre_charac_cub}
and Theorem~\ref{preacc}.  \epf

\begin{cor} 
  There exists a left determined model structure with respect to the
  class of cofibrations between cubical transition systems. The
  adjunction $\cts \leftrightarrows \whdts$ is a Quillen adjunction.
  All objects of $\cts$ are cofibrant.
\end{cor}

\bpf The class of cofibrations between cubical transition systems is
generated by a set $\I^{\cts}$ by Theorem~\ref{catlem}. The segment $V$
is cubical by Proposition~\ref{charac_cub}. The other hypotheses of
Theorem~\ref{constr_submodel} are easy to check. Hence the proof is
complete.  \epf

Proposition~\ref{pre_charac_cub} has a consequence which will not be
used in the paper but which is worth mentioning anyway. This is about
an explicit description of the coreflector from $\whdts$ to
$\cts$. 

\bd Let $X$ be a weak HDTS. A $(n+1)$-transition
$(\alpha,u_1,\dots,u_{n+1},\beta)$ of $X$ is {\rm divisible} if either
$n=0$ or there exists a state $\gamma$ such that the tuples
$(\alpha,u_1,\dots,u_{p},\gamma)$ and
$(\gamma,u_{p+1},\dots,u_{n+1},\beta)$ are two divisible transitions
of $X$ for some $p \geq 1$. \ed

\bp \label{coreflector-explicit} Let $X$ be a weak HDTS. The image
$\overline{X}$ of $X$ by the coreflector is the weak HDTS having the
same states as $X$, having as set of actions the actions of $X$ which
are used in a $1$-transition (in the sense of Lemma~\ref{used}) and
having as set of transitions the divisible transitions.  \ep

\bpf It is clear by Proposition~\ref{pre_charac_cub} that all
transitions of $\overline{X}$ are divisible. Conversely, let
$(\alpha,u_1,\dots,u_n,\beta)$ be a divisible transition of $X$. Then
the corresponding map \[C_n[\mu(u_1),\dots,\mu(u_n)]^{ext}
\longrightarrow X\] factors as a
composite \[C_n[\mu(u_1),\dots,\mu(u_n)]^{ext} \longrightarrow
C_n[\mu(u_1),\dots,\mu(u_n)] \longrightarrow X.\] Therefore every
divisible transition belongs to a subcube.  \epf

\section{First Cattani-Sassone axiom and weakly equivalent cubical transition systems}
\label{weakisiso}

From now on, we work in the category of cubical transition systems
$\cts$. So $\cof = \cof_{\cts}$, $\inj = \inj_{\cts}$, $\cell =
\cell_{\cts}$. The localizer (with respect to the class of
cofibrations of cubical transition systems) generated by a set
$\mathcal{S}$ is denoted by $\W(\mathcal{S})$.

We want to characterize the weak equivalences of the left determined
model structure of cubical transition systems. The following axiom,
introduced in \cite{hdts}, will be useful.

\bd \label{csa1} A cubical transition system satisfies the {\rm First
  Cattani-Sassone axiom} (CSA1) if for every transition
$(\alpha,u,\beta)$ and $(\alpha,u',\beta)$ such that the actions $u$
and $u'$ have the same label in $\Sigma$, one has $u = u'$. \ed

The axiom CSA1 used by Cattani and Sassone in their paper
\cite{MR1461821} is even stronger, but we do not need this stronger
form. In our language, their stronger form states that if
$(\alpha,u_1,\dots,u_n,\beta)$ and $(\alpha,u'_1,\dots,u'_n,\beta)$
are two $n$-dimensional transitions with $\mu(u_i) = \mu(u'_i)$ for
$1\leq i\leq n$, then one has $(\alpha,u_1,\dots,u_n,\beta) =
(\alpha,u'_1,\dots,u'_n,\beta)$.

\bp The full subcategory of cubical transition systems satisfying CSA1
is a full reflective subcategory of $\cts$. \ep

\bpf The category of cubical transition systems satisfying CSA1 is a
small-orthogonality class of $\cts$. Indeed a cubical transition
system satisfies CSA1 if and only if it is orthogonal to the set of
maps $C_1[x] \sqcup_{\{0_1,1_1\}} C_1[x] \longrightarrow C_1[x]$ for
$x$ running over $\Sigma$. The proof goes exactly as in
\cite[Corollary~5.7]{hdts}. \epf

\begin{nota} Let us denote by $\CSA_1$  the reflector. \end{nota}

\bp \label{homotopy_is_equality} Let $Y$ be a cubical transition
system satisfying CSA1.  Let $X$ be a cubical transition system. Then
two homotopy equivalent maps $f,g:X\rightarrow Y$ are equal. In other
terms, each of the two canonical maps $X \rightarrow X\p V$ induces a
bijection $\cts(X\p V,Y) \iso \cts(X,Y)$. \ep

\bpf The cubical transition system $X\p V$ is calculated in the proof
of Proposition~\ref{exp}. Let us recall the results. The cubical
transition system $X\p V$ and $X$ have the same states. If $L$ is the
set of actions of $X$, then $L\p \{0,1\}$ is the set of actions of
$X\p V$ and the labelling map is the composite $L\p \{0,1\}
\rightarrow L \rightarrow \Sigma$. Finally, a tuple
$(\alpha,(u_1,\epsilon_1),\dots,(u_n,\epsilon_n),\beta)$ for
$\epsilon_1,\dots,\epsilon_n \in \{0,1\}$ is a transition of $X\p V$
if and only if the tuple $(\alpha,u_1,\dots,u_n,\beta)$ is a
transition of $X$.

Let us consider a homotopy $H : X \p V \rightarrow Y$ between two maps
$f$ and $g$ from $X$ to $Y$. Since $X\p V$ and $X$ have the same
states, $f_0 = g_0 = H_0$, i.e. $f$ and $g$ coincide on states. Let
$u$ be an action of $X$. Since $X$ is injective with respect to the
map $\underline{\mu(u)} \longrightarrow C_1[\mu(u)]$ by
Theorem~\ref{preacc}, there exists a transition $(\alpha,u,\beta)$ of
$X$. So the tuples $(\alpha,(u,0),\beta)$ and $(\alpha,(u,1),\beta)$
are two transitions of $X\p V$. Therefore
$(H_0(\alpha),\widetilde{H}(u,0),H_0(\beta))$ and
$(H_0(\alpha),\widetilde{H}(u,1),H_0(\beta))$ are two transitions of
$Y$. By CSA1, one has $\widetilde{f}(u) = \widetilde{H}(u,0) =
\widetilde{H}(u,1) = \widetilde{g}(u)$. Hence $f = g$.  \epf

\begin{cor} \label{CSA1_PathObject} Let $T$ be a cubical transition
  system satisfying CSA1.  Then there is the canonical isomorphism
  $T^V \iso T$ in $\cts$~\footnote{The weak HDTS $(T^V)^{\whdts}$ (the
    right adjoint being calculated in $\whdts$) is not isomorphic to
    $T$; the calculations in the proof of Proposition~\ref{exp} show
    that the two weak HDTS have a different set of actions,
    $L\p_\Sigma L$ for $(T^V)^{\whdts}$ if $L$ is the set of actions
    of $T$.} \end{cor}

\bp \label{critere} Let $T$ be a cubical transition system such that
$T^V \iso T$ (in $\cts$). Then one has:
\begin{enumerate}
\item $T$ is orthogonal to every map of the form $f\star
  \gamma^\epsilon$ with $\epsilon=0,1$ and with $f$ any map of cubical
  transition systems.
\item $T$ is injective with respect to a map of the form $f\star
  \gamma$ with $f$ a map of cubical transition systems if and only if
  for every diagram of the form
\[
\xymatrix{
X \ar@{->}[r]^-{g}\ar@{->}[d]^-{f}  & T\\
Y \ar@{-->}[ru]_-{k}&}
\] there exists at most one lift $k$.
\item $T$ is injective with respect to every map of the form $f\star
  \gamma$ with $f$ a map of cubical transition systems such that $f_0$
  and $\widetilde{f}$ are onto~\footnote{In fact, this assertion holds
    whenever $f$ is an epimorphism.}.
\item $T$ is injective with respect to every map of the form $(f\star
  \gamma) \star \gamma$ where $f$ is a map of cubical transition
  systems.
\end{enumerate}
\ep 

\bpf By adjunction, $T$ is injective with respect to a map of the form
$f\star \gamma^\epsilon$ if and only if $f$ satisfies the LLP with
respect to the map $\pi_\epsilon : T^V \rightarrow T$ which is an
isomorphism. Hence the first assertion. 

By adjunction again, $T$ is injective with respect to a map of the
form $f\star \gamma$ if and only if $f$ satisfies the LLP with respect
to the canonical map $\pi: T^V \rightarrow T\p T$ which turns out to
be the diagonal.  Two lifts $k_1$ and $k_2$ in the diagram 
\[
\xymatrix{
X \ar@{->}[r]^-{g}\ar@{->}[d]^-{f}  & T\\
Y \ar@{-->}[ru]_-{k_1,k_2}&}
\]
give rise to the commutative diagram of solid arrows
\[
\xymatrix{
  X \fR{g} \fD{f} && T^V \fD{} \\
  && \\
  Y \fR{(k_1,k_2)} \ar@{-->}[rruu]^-{k}&& T\p T.}
\]  
One deduces $k_1 = k = k_2$. Conversely, let us suppose that there is
always at most one lift $k$ in the diagram 
\[
\xymatrix{
X \ar@{->}[r]^-{g}\ar@{->}[d]^-{f}  & T\\
Y \ar@{-->}[ru]_-{k}&}
\]
Consider a commutative diagram of solid arrows of the form
\[
\xymatrix{
  X \fR{g} \fD{f} && T^V \iso T \fD{} \\
  && \\
  Y \fR{(k_1,k_2)} && T\p T.}
\] 
Then $k_1 = k_2$ and therefore $T$ is $(f\star
\gamma)$-injective. Hence the second assertion.

Let us suppose now that $f$ is a map of cubical transition systems
such that $f_0$ and $\widetilde{f}$ are onto. Let $k_1$ and $k_2$ be
two lifts. Then $\omega(k_1)\omega(f) = \omega(g) =
\omega(k_2)\omega(f)$. So $\omega(k_1) = \omega(k_2)$. Since the
forgetful functor $\omega$ is faithful, one deduces that $k_1 =
k_2$. Hence the third assertion.

Let $f : X \rightarrow X'$ be a map of cubical transition systems with
$X = (S,\mu : L \rightarrow \Sigma, T)$ and $X' = (S',\mu' : L'
\rightarrow \Sigma, T')$. The map $f \star \gamma$ is obtained by
considering the commutative diagram of solid arrows
\[
\xymatrix{
  X \sqcup X \fD{} \fR{} && \cyl(X) \fD{} \\
  && \\
  X' \sqcup X' \fR{} && \cyl(X')}
\] and by using the universal property of the pushout, giving the map
\[f \star \gamma : (X'\sqcup X') \sqcup_{X\sqcup X} \cyl(X)
\longrightarrow \cyl(X')\] The latter map induces on the set of states
the map $(S' \sqcup S') \sqcup_{S \sqcup S} S \iso S' \sqcup_S S'
\rightarrow S'$ which is onto, and on the set of actions the map $(L'
\sqcup L') \sqcup_{L\sqcup L} (L\sqcup L) \iso L' \sqcup L'
\rightarrow L' \sqcup L'$ which is onto as well. So the fourth
assertion is a consequence of the third one. \epf

\bp \label{fibrant0} Let $\mathcal{S}$ be a set of maps of cubical
transition systems.  Let $T$ be a cubical transition system satisfying
CSA1. Then $T$ is $\Lambda(\cyl,\mathcal{S},\I^{\cts})$-injective if
and only if $T$ is $\mathcal{S}$-orthogonal.  \ep

\bpf If $T$ is $\Lambda(\cyl,\mathcal{S},\I^{\cts})$-injective, then
it is $\Lambda^0(\cyl,\mathcal{S},\I^{\cts})$-injective, and therefore
$\mathcal{S}$-injective. Such a $T$ is also
$\Lambda^1(\cyl,\mathcal{S},\I^{\cts})$-injective with
\[\Lambda^1(\cyl,\mathcal{S},\I^{\cts}) =
\Lambda^0(\cyl,\mathcal{S},\I^{\cts}) \star \gamma.\] Therefore $T$ is
$\mathcal{S}$-orthogonal by Proposition~\ref{critere} (2). Conversely,
let us suppose that $T$ is $\mathcal{S}$-orthogonal. Then $T$ is
$\Lambda^0(\cyl,\mathcal{S},\I^{\cts})$-injective by
Proposition~\ref{critere} (1). By Proposition~\ref{critere} (2) and
(1), $T$ is $\Lambda^1(\cyl,\mathcal{S},\I^{\cts})$-injective as
well. The injectivity with respect to
$\Lambda^n(\cyl,\mathcal{S},\I^{\cts})$ for $n\geq 2$ is a consequence
of Proposition~\ref{critere} (4).\epf

Hence the theorems: 

\bp \label{fibrant} Every cubical transition system satisfying CSA1 is
fibrant in the left determined model structure of $\cts$. \ep

\bpf The statement is a corollary of Proposition~\ref{fibrant0} with
$\mathcal{S} = \varnothing$. \epf

\bp \label{weak_is_equal} Two cubical transition systems satisfying
CSA1 are weakly equivalent if and only if they are isomorphic.  \ep

\bpf Let $f: X \rightarrow Y$ be a weak equivalence between two
cubical transition systems satisfying CSA1. Since $X$ and $Y$ are both
cofibrant and fibrant by Proposition~\ref{fibrant}, there exists a map
$g:Y \rightarrow X$ such that $f\circ g$ is homotopy equivalent to
$\id_Y$ and such that $g\circ f$ is homotopy equivalent to $\id_X$.
So by Proposition~\ref{homotopy_is_equality}, $f\circ g = \id_Y$ and
$g\circ f = \id_X$. Hence $X$ and $Y$ are isomorphic.
\epf

\bth \label{weak-equivalence-known} The reflector $\CSA_1$ detects the
weak equivalences of the left determined model structure of $\cts$. In
other terms, a map $f$ of cubical transition systems is a weak
equivalence in the left determined model structure of $\cts$ if and
only if $\CSA_1(f)$ is an isomorphism. \eth

In particular, this theorem means that two cubical transition systems
interpreting two process names in a process algebra are weakly
equivalent in this model structure if and only if they are
isomorphic. See \cite{hdts} for further details.

\bpf By Proposition~\ref{weak_is_equal}, it suffices to prove that for
every cubical transition system $X$, the unit $X \rightarrow
\CSA_1(X)$ is a weak equivalence in the left determined model
structure of $\cts$. An object $X$ is orthogonal to a map of the form
$C_1[x] \sqcup_{\{0_1,1_1\}} C_1[x] \longrightarrow C_1[x]$ for $x \in
\Sigma$ if and only if it is injective with respect to it since this
map is an epimorphism. So the map $X \rightarrow \CSA_1(X)$ is
obtained by factoring the canonical map $X \rightarrow \mathbf{1}$
(from $X$ to the terminal object) as a composite $X \rightarrow
\CSA_1(X) \rightarrow \mathbf{1}$ where the left-hand map belongs to
$\cell_{\cts}(\mathcal{U})$ and the right-hand map belongs to
$\inj_{\cts}(\mathcal{U})$ where \[\mathcal{U} = \{C_1[x]
\sqcup_{\{0_1,1_1\}} C_1[x] \longrightarrow C_1[x] \mid x\in
\Sigma\}.\] So it suffices to prove that every pushout of a map of the
form $C_1[x] \sqcup_{\{0_1,1_1\}} C_1[x] \rightarrow C_1[x]$ for $x\in
\Sigma$ is a weak equivalence of the left determined model structure
of $\cts$. The identity of $C_1[x]$ factors as a composite \[C_1
\longrightarrow C_1[x] \sqcup_{\{0_1,1_1\}} C_1[x] \longrightarrow
C_1[x].\] By the calculation made in the proof of
Proposition~\ref{exp}, there is the isomorphism $C_1[x] \p V \iso
C_1[x] \sqcup_{\{0_1,1_1\}} C_1[x]$. Hence the left-hand map is a weak
equivalence, and also the right-hand map by the two-out-of-three
axiom. Consider a pushout diagram of the form
\[
\xymatrix{
C_1[x] \sqcup_{\{0_1,1_1\}} C_1[x] \fR{\phi} \fD{} && X \fD{f} \\
&& \\
C_1[x] \fR{} && \cocartesien Y}
\]
The cubical transition system $C_1[x] \sqcup_{\{0_1,1_1\}} C_1[x]$
contains two actions $x_1$ and $x_2$ labelled by $x$.  There are two
mutually exclusive cases. Either $\widetilde{\phi}(x_1) =
\widetilde{\phi}(x_2)$ or $\widetilde{\phi}(x_1) \neq
\widetilde{\phi}(x_2)$. In the first case, the commutative square above factors as a composite 
of commutative squares 
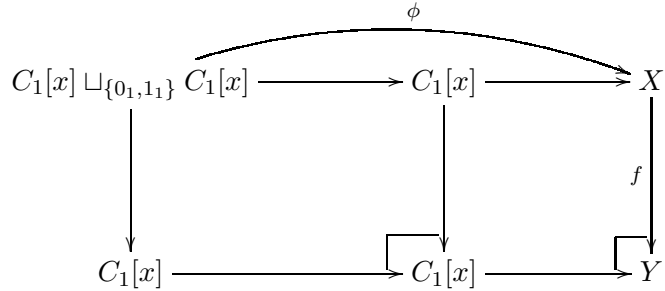
\[
\xymatrix{
C_1[x] \sqcup_{\{0_1,1_1\}} C_1[x] \fR{}\ar@/^20pt/[rrrr]^-{\phi} \fD{} && C_1[x]\fR{} \fD{} &&\fD{f} X\\
&& &&\\
C_1[x] \fR{} && \cocartesien C_1[x]\fR{} && \cocartesien Y}
\]
Hence $X \iso Y$. In the second case, $\phi$ is one-to-one on actions,
i.e. a cofibration of cubical transition systems. In that case, $f$ is
a weak equivalence since the left determined model structure of $\cts$
is left proper. So the map $X \rightarrow \CSA_1(X)$ is a transfinite
composition of weak equivalences. The class of weak equivalences of a
combinatorial model category is always accessible accessibly-embedded
by e.g. \cite[Corollary~A.2.6.6]{Lurie}. Hence a transfinite
composition of weak equivalences is always a weak equivalence. The
proof is complete.  \epf

\begin{cor} \label{not-enough} The counit map $p_x : \cub(\dd{x})
  \longrightarrow \dd{x}$ is not a weak equivalence in the left
  determined model structure of $\cts$.
\end{cor}

Corollary~\ref{not-enough} shows that this model structure is really
minimal. Even cubical transition systems having the same cubes may be
not weakly equivalent.  The next section explains how it is possible
to add weak equivalences so that two cubical transition systems
containing the same cubes after simplification of the labelling are
always weakly equivalent.

\section{Bousfield localization with respect to the cubification functor}
\label{loc}

Let us denote by $\W_{\cub}$ the smallest localizer generated by the
class of maps of cubical transition systems $f:X\rightarrow Y$ such
that $\cub(f)$ is a weak equivalence in the left determined model
structure of $\cts$.  We want to prove that it is small, more
precisely that it is generated by the set of maps $\mathcal{S} =
\{p_x:C_1[x] \sqcup C_1[x] \rightarrow \dd{x}\mid x\in \Sigma\}$.

Let us prove first that the two functors $\cub(-)$ and $\CSA_1(-)$
commute with one another.

\bp \label{cub-preserves-CSA1} Let $X$ be a cubical transition
system. Then there exists a natural isomorphism $\CSA_1(\cub(X)) \iso
\cub(\CSA_1(X))$. \ep

\bpf That $\CSA_1(X)$ satisfies CSA1 means that for every $x\in
\Sigma$, the map $C_1[x] \sqcup_{\{0_1,1_1\}} C_1[x] \rightarrow
C_1[x]$ induces a bijection 
\[\cts(C_1[x],\CSA_1(X)) \iso \cts(C_1[x]
\sqcup_{\{0_1,1_1\}} C_1[x],\CSA_1(X)).\] By
Proposition~\ref{colim_gen}, the functor $\cub$ is right adjoint to
the inclusion functor of the full subcategory of $\cts$ generated by
the cubes $C_n[x_1,\dots,x_n]$ for $n\geq 0$ and $x_1,\dots,x_n
\in\Sigma$ into $\cts$. Both $C_1[x]$ and $C_1[x] \sqcup_{\{0_1,1_1\}}
C_1[x]$ are colimits of cubes. Therefore one has the bijections
\begin{align*}
  &\cts(C_1[x] \sqcup_{\{0_1,1_1\}}
  C_1[x],\cub(\CSA_1(X))) & \\
  &  \iso \cts(C_1[x]
  \sqcup_{\{0_1,1_1\}} C_1[x],\CSA_1(X)) & \hbox{ by adjunction} \\
  & \iso \cts(C_1[x],\CSA_1(X)) & \hbox{ since $\CSA_1(X)$ satisfies CSA1}\\
  & \iso
  \cts(C_1[x],\cub(\CSA_1(X))) & \hbox{ by adjunction again.}
\end{align*}
Hence $\cub(\CSA_1(X))$ satisfies CSA1.  Therefore the canonical map
\[\xymatrix{\cub(X) \fR{\cub(\phi_X)} &&  \cub(\CSA_1(X))}\] factors
uniquely as a composite
\[\xymatrix{\cub(X) \fR{\phi_{\cub(X)}}&& \CSA_1(\cub(X)) \fR{\psi_X}&& 
\cub(\CSA_1(X)).}\] The functors $\cub$ and $\CSA_1$ preserve states. So
the map $\psi_X$ is a bijection on states. The map $\psi_X$ is also
surjective on actions and on transitions since any of them comes
respectively from an action or a transition of $\cub(X)$.

It remains to understand why the map $\psi_X$ is one-to-one on actions
for the proof to be complete. Consider the commutative diagram of
cubical transition systems of Figure~\ref{square1}.  Since the cubical
transition systems of the bottom line of Figure~\ref{square1} satisfy
CSA1, this square factors uniquely as a composite of commutative
squares as in Figure~\ref{square2}. Let $u_1$ and $u_2$ be two actions
of $\CSA_1(\cub(X))$ such that ${\psi_X}(u_1) = {\psi_X}(u_2) = u$. Let
$u'_1$ and $u'_2$ be two actions of $\cub(X)$ such that
${\phi_{\cub(X)}}(u'_1) = u_1$ and ${\phi_{\cub(X)}}(u'_2) = u_2$. Let
$v'_1 = p_X(u'_1)$, $v'_2 = p_X(u'_2)$, $v_1 = \CSA_1(p_X)(u_1)$, $v_2
= \CSA_1(p_X)(u_2)$ and finally $v = p_{\CSA_1(X)}(u)$\footnote{We
  denote in the same way a map of cubical transition systems $f$ and
  the set map $\widetilde{f}$ between actions in order to not overload
  the notations.}  By commutativity of the diagram, we obtain $v_1 =
v_2 = v$. By construction of the functor $\CSA_1(-)$, there exist two
states $\alpha$ and $\beta$ such that the triple $(\alpha,v'_1,\beta)$
and $(\alpha,v'_2,\beta)$ are two transitions of $X$. Therefore by
definition of $\cub$, the two triples $(\alpha,u'_1,\beta)$ and
$(\alpha,u'_2,\beta)$ are two transitions of $\cub(X)$. So $u_1 =
u_2$ since $\CSA_1(\cub(X))$ satisfies CSA1.  \epf

\begin{figure}
\[
\xymatrix
{
\cub(X) \fR{p_X}\fD{\cub(\phi_X)} && X \fD{\phi_X} \\
&&\\
\cub(\CSA_1(X)) \fR{p_{\CSA_1(X)}} && \CSA_1(X)
}
\]
\caption{Composition of $\cub$ and $\CSA_1$ (I)}
\label{square1}
\end{figure}

\begin{figure}
\[
\xymatrix
{
\cub(X) \fR{p_X}\fD{\phi_{\cub(X)}} && X \fD{\phi_X} \\
&&\\
\CSA_1(\cub(X)) \fR{\CSA_1(p_X)} \fD{\psi_X}&&  \CSA_1(X) \ar@{=}[dd]\\
&&\\
\cub(\CSA_1(X)) \fR{p_{\CSA_1(X)}} && \CSA_1(X)
}
\]
\caption{Composition of $\cub$ and $\CSA_1$ (II)}
\label{square2}
\end{figure}

\bp The functor $\cub : \cts \rightarrow \cts$ preserves weak equivalences. \ep

\bpf Let $f$ be a weak equivalence of $\cts$. Then $\CSA_1(f)$ is an
isomorphism by Theorem~\ref{weak-equivalence-known}. So
$\CSA_1(\cub(f))$ is an isomorphism by
Proposition~\ref{cub-preserves-CSA1}. Therefore by
Theorem~\ref{weak-equivalence-known} again, $\cub(f)$ is a weak
equivalence of $\cts$.  \epf

\begin{cor} Every weak equivalence of $\cts$ belongs to
  $\W_{\cub}$. \end{cor}

\bp \label{ll1} Let $X$ be a cubical transition system. The counit
$p_X : \cub(X) \rightarrow X$ is a transfinite composition of pushouts
of the maps $p_x:C_1[x] \sqcup C_1[x] \rightarrow \dd{x}$ for $x$
running over $\Sigma$.  \ep

\bpf We already know that the map $p_X : \cub(X) \rightarrow X$ is
bijective on states. let $u$ be an action of $X$. Since $X$ is
cubical, there exists a $1$-transition $(\alpha,u,\beta)$ of $X$,
which corresponds to a map $C_1[\mu(u)] \rightarrow X$. Hence the map
$p_X : \cub(X) \rightarrow X$ is onto on actions. Let
$(\alpha,u_1,\dots,u_n,\beta)$ be a transition of $X$, which
corresponds to a map $C_n[\mu(u_1),\dots,\mu(u_n)]^{ext} \rightarrow
X$. Since $X$ is cubical, the latter map factors as a composite
$C_n[\mu(u_1),\dots,\mu(u_n)]^{ext} \rightarrow
C_n[\mu(u_1),\dots,\mu(u_n)] \rightarrow X$ by
Theorem~\ref{preacc}. Hence the map $p_X : \cub(X) \rightarrow X$ is
onto on transitions. Let us factor the map $p_X$ as a composite
$\cub(X) \rightarrow Z \rightarrow X$ where the left-hand map belongs
to $\cell(\mathcal{S})$ and the right-hand map belongs to
$\inj(\mathcal{S})$. The right-hand map $ g: Z\rightarrow X$ is still
bijective on states, and onto on actions and transitions. Let $u_1$
and $u_2$ be two actions of $Z$ mapped to the same action $u$ of
$X$. Then $\mu(u_1) = \mu(u_2) = \mu(u) = x$.  Let us suppose that the
action $u_1$ is used in a transition $(\alpha_1,u_1,\beta_1)$, and the
action $u_2$ in a transition $(\alpha_2,u_2,\beta_2)$ of $Z$.  Then
consider the commutative diagram of cubical transition systems
\[
\xymatrix
{
C_1[x] \sqcup C_1[x] \fR{} \fD{} && Z \fD{}\\
&& \\
\dd{x} \ar@{-->}[rruu]^-{\ell}\fR{} && X,
}
\]
where each copy of $C_1[x]$ corresponds to one of the two transitions
$(\alpha_i,u_i,\beta_i)$. The existence of the lift $\ell$ implies
that $u_1 = u_2$. So the map $ g: Z\rightarrow X$ is one-to-one on
actions.  Finally, let $(\alpha,u_1,\dots,u_n,\beta)$ and
$(\alpha',u'_1,\dots,u'_n,\beta')$ be two transitions of $Z$ mapped to
the same transition of $X$. Then $\alpha = \alpha'$, $\beta = \beta'$
and $u_i = u'_i$ for $1\leq i \leq n$ since $ g: Z\rightarrow X$ is
bijective on states and actions. So $g$ is one-to-one on
transitions. Therefore $g$ is an isomorphism.  \epf

\bp \label{ll2} Every map of $\cell(\mathcal{S})$ belongs to the
localizer generated by $\mathcal{S}$, i.e.  $\cell(\mathcal{S})
\subset \W(\mathcal{S})$. \ep

\bpf Note that $p_x : C_1[x] \sqcup C_1[x] \rightarrow \dd{x}$ is not
a cofibration so we cannot use the fact that the class of trivial
cofibrations is closed under pushout and transfinite compositions.  By
Theorem~\ref{constr_submodel}, the class of maps $\W(\mathcal{S})$ is
the class of weak equivalences of a model structure on $\cts$. It is
actually the class of weak equivalences of the Bousfield localization
of the left determined model structure of $\cts$ by
$\mathcal{S}$. Since all objects are cofibrant, it is left
proper. Consider a pushout diagram of the form
\[
\xymatrix{
C_1[x] \sqcup C_1[x] \fR{\phi} \fD{p_x} && X \fD{}\\
&& \\
\dd{x} \fR{} && Y. \cocartesien
}
\]
There are two mutually exclusive cases. The map $\phi$ takes the two
actions of $C_1[x] \sqcup C_1[x]$ to two different actions. Then
$\phi$ is a cofibration and $X\rightarrow Y$ belongs to
$\W(\mathcal{S})$ by left properness. Or $\phi$ takes the two actions
of $C_1[x] \sqcup C_1[x]$ to the same action. Then $X\iso Y$ (the
argument is similar to the one used in the proof of
Theorem~\ref{weak-equivalence-known}). So in the two cases, the
right-hand vertical map belongs to $\W(\mathcal{S})$.  The proof is
complete by \cite[Corollary~A.2.6.6]{Lurie} since $\W(\mathcal{S})$ is
closed under transfinite composition.  \epf

Hence the theorem: 

\bth \label{ll3} One has the equality of localizers $\W_{\cub} =
\W(\mathcal{S})$. \eth

\bpf The map $\cub(p_x)$ is an isomorphism by
Proposition~\ref{facile_1} and Proposition~\ref{colim_gen}. So
$\mathcal{S} \subset \W_{\cub}$. Hence the first inclusion
$\W(\mathcal{S}) \subset \W_{\cub}$. Let $f : X\rightarrow Y$ be a map
of cubical transition systems such that $\cub(f)$ is a weak
equivalence of the left determined model structure of $\cts$.
Consider the commutative diagram
\[
\xymatrix{
\cub(X) \fR{\cub(f)} \fD{} && \cub(Y) \fD{} \\
&& \\
X \fR{f} && Y.}
\]
The vertical maps belong to $\W(\mathcal{S})$ by Proposition~\ref{ll1}
and Proposition~\ref{ll2}. By hypothesis, the top horizontal map is a
weak equivalence of $\cts$, and therefore belongs to $\W(\mathcal{S})$
as well. Hence by the two-out-of-three property, $f : X\rightarrow Y$
belongs to $\W(\mathcal{S})$. We obtain the second inclusion
$\W_{\cub} \subset \W(\mathcal{S})$. \epf

\begin{cor} \label{localization} The Bousfield localization of the
  left determined model structure of $\cts$ with respect to the
  functor $\cub$ exists. \end{cor}

The weak factorization system $(\cof(\mathcal{S}),\inj(\mathcal{S}))$
gives rise to a functor $\bl_{\mathcal{S}} : \cts \rightarrow
\cts$. It is defined by functorially factoring the map $X \rightarrow
\mathbf{1}$ as a composite $X \rightarrow \bl_{\mathcal{S}}(X)
\rightarrow \mathbf{1}$ where the left-hand map belongs to
$\cell(\mathcal{S})$ and the right-hand map belongs to
$\inj(\mathcal{S})$.

\begin{rem} The labelling map is one-to-one for every cubical
  transition system of the form $\bl_{\mathcal{S}}(X)$. \end{rem}

A remarkable consequence of this fact is that for every map $f:X
\rightarrow Y$ of cubical transition systems, the map
$\bl_{\mathcal{S}}(f) : \bl_{\mathcal{S}}(X) \rightarrow
\bl_{\mathcal{S}}(Y)$ is a cofibration. So $\bl_{\mathcal{S}}(f)$ is a
cofibrant replacement of $f$ in $\bl_{\mathcal{S}}(\cts)$ by
Proposition~\ref{ll2} and Theorem~\ref{ll3}.

\bp \label{injS-CSA1} Every cubical transition system of
$\inj(\mathcal{S})$ satisfies CSA1. \ep

\bpf Let $(\alpha_1,u_1,\beta_1)$ and $(\alpha_2,u_2,\beta_1)$ be two
$1$-transitions of a cubical transition system injective with respect
to $\mathcal{S}$ with $\mu(u_1) = \mu(u_2)$. Then $u_1 = u_2$, and
this is still true if $\alpha_1 = \alpha_2$ and $\beta_1 =
\beta_2$. Hence CSA1 is satisfied.  \epf

The weak equivalences of this Bousfield localization have a nice
characterization.

\bth \label{detect} A map of cubical transition systems $f : X
\rightarrow Y$ belongs to $\W(\mathcal{S})$ if and only if
$\bl_{\mathcal{S}}(f) : \bl_{\mathcal{S}}(X) \iso
\bl_{\mathcal{S}}(Y)$ is an isomorphism. In other terms, the functor
$\bl_{\mathcal{S}}$ detects the weak equivalences of this Bousfield
localization.  \eth

\bpf If $\bl_{\mathcal{S}}(f)$ is an isomorphism, $f$ belongs to
$\W(\mathcal{S})$ by Proposition~\ref{ll2} and by the two-out-of-three
property.  Conversely, suppose that $f \in \W(\mathcal{S})$. By
Proposition~\ref{injS-CSA1} and Corollary~\ref{CSA1_PathObject}, one
has $\bl_{\mathcal{S}}(X)^V = \bl_{\mathcal{S}}(X)$ and
$\bl_{\mathcal{S}}(Y)^V = \bl_{\mathcal{S}}(Y)$. The maps of
$\mathcal{S}$ are onto on states and actions.  So by
Proposition~\ref{fibrant0}, $\bl_{\mathcal{S}}(X)$ and
$\bl_{\mathcal{S}}(Y)$ are fibrant in the Bousfield localization since
if they are orthogonal to the maps of $\mathcal{S}$. By the
two-out-of-three property, $\bl_{\mathcal{S}}(f)$ is therefore a weak
equivalence between two cofibrant-fibrant objects in the Bousfield
localization $\bl_{\mathcal{S}}(\cts)$ of the left determined model
structure of $\cts$ by the maps of $\mathcal{S}$. By
\cite[Theorem~3.2.13]{ref_model2}, the map $\bl_{\mathcal{S}}(f)$ is
then a weak equivalence of the left determined model structure of
$\cts$. Since $\bl_{\mathcal{S}}(X)$ and $\bl_{\mathcal{S}}(Y)$
satisfy CSA1 by Proposition~\ref{injS-CSA1}, the map
$\bl_{\mathcal{S}}(f)$ is an isomorphism by
Theorem~\ref{weak-equivalence-known}.  \epf

So in the Bousfield localization $\bl_{\mathcal{S}}(\cts)$, two
cubical transition systems are weakly equivalent if they have the same
cubes after simplification of the labelling. It is actually possible
to prove better:

\bth We have:
\begin{enumerate}
\item The functor $\bl_{\mathcal{S}} : \cts \rightarrow \cts$
  induces a functor from $\cts$ to the full reflective subcategory
  $\mathcal{S}^\perp$ of cubical transition systems consisting of
  $\mathcal{S}$-orthogonal objects.
\item For every $\mathcal{S}$-orthogonal cubical transition system
  $Y$, there is a natural isomorphism $Y \iso \bl_{\mathcal{S}}(Y)$.
\item The functor $\bl_{\mathcal{S}}$ is left adjoint to the inclusion
  functor $\mathcal{S}^\perp \subset \cts$.
\item Every map between $\mathcal{S}$-orthogonal cubical transition
  systems is a cofibration of cubical transition systems. Every
  $\mathcal{S}$-orthogonal cubical transition system is cofibrant and
  fibrant in $\bl_{\mathcal{S}}(\cts)$.
\item The homotopy category of $\bl_{\mathcal{S}}(\cts)$ is equivalent
  to $\mathcal{S}^\perp$.
\end{enumerate}
\eth

\bpf (1) comes from the definition of $\bl_{\mathcal{S}}$ and from the
fact that $\mathcal{S}$-injective is equivalent to
$\mathcal{S}$-orthogonal since every map of $\mathcal{S}$ is an
epimorphism.  One has a natural isomorphism $\bl_{\mathcal{S}}(Y) \iso
Y$ for every $\mathcal{S}^\perp$-orthogonal cubical transition system
$Y$ since every pushout $Y \rightarrow Z$ of a map of the form
$p_x:C_1[x] \sqcup C_1[x] \rightarrow \dd{x}$ for $x\in \Sigma$ is an
isomorphism, hence (2).  For every $\mathcal{S}^\perp$-orthogonal
cubical transition system $Y$, the canonical map $Y \rightarrow
\mathbf{1}$ satisfies the RLP with respect to every map of
$\cell(\mathcal{S})$, in particular with respect to every map $X
\rightarrow \bl_{\mathcal{S}}(X)$ for every cubical transition system
$X$. Moreover, every map of $\cell(\mathcal{S})$ is bijective on
states and onto on actions; so every map of $\cell(\mathcal{S})$ is an
epimorphism. So $\cell(\mathcal{S})$-injective is equivalent to
$\cell(\mathcal{S})$-orthogonal. This means that every map $X
\rightarrow Y$ from a cubical transition system $X$ to an
$\mathcal{S}$-orthogonal cubical transition system $Y$ factors
uniquely as a composite $X \rightarrow \bl_{\mathcal{S}}(X)\rightarrow
Y$, hence (3). (4) is explained in the proof of Theorem~\ref{detect}.
The functor $\bl_{\mathcal{S}}: \cts \rightarrow \cts$ factors
uniquely as a composite $\cts \rightarrow \bl_{\mathcal{S}}(\cts)
\rightarrow \mathcal{S}^\perp$ by Theorem~\ref{detect} and by the
universal property of the categorical localization. There is a natural
isomorphism $X \rightarrow \bl_{\mathcal{S}}(X)$ in
$\bl_{\mathcal{S}}(\cts)$ by Proposition~\ref{ll2} for every object of
$\cts$. And there is a natural isomorphism $Y\iso
\bl_{\mathcal{S}}(Y)$ for every $\mathcal{S}$-orthogonal object since
$\mathcal{S}$-injective is equivalent to
$\mathcal{S}$-orthogonal. Hence (5).  \epf

\section{Weak equivalence and  bisimulation}
\label{bisim}

This last section sketches the link between these homotopical
constructions and bisimulation. Let us introduce bisimulations with
open maps as in \cite{0856.68067}. The link between bisimulation and
homotopy will be the subject of future works. Indeed, the definition
of open maps taken here is very restrictive since a good definition
requires a more general notion of paths (cf. \cite{fahrenberg05-hda}
for further explanations). The purpose of this section is only to have
an idea of what it is possible to do with these homotopical
constructions.

Let $\mathcal{P}$ be a subset of the set of cubes
$\{C_n[x_1,\dots,x_n]\mid n\geq 0,x_1,\dots,x_n\in\Sigma\}$. The
elements of $\mathcal{P}$ are called \emph{calculation paths}.

\bd A map $f : X \rightarrow Y$ is {\rm $\mathcal{P}$-open} if every
commutative square of solid arrows
\[
\xymatrix{
\{0_n\} \fR{} \fD{} && X \fD{} \\
&& \\
P \fR{} \ar@{-->}[rruu]^-{k}&& Y}
\]
as a lift $k$ for every $P\in \mathcal{P}$, i.e. $f$ satisfies the RLP
with respect to the inclusion $\{0_n\} \subset P$.
\ed

\bd Two cubical transition systems $X$ and $Y$ are {\rm
  $\mathcal{P}$-bisimilar} if there exists a cubical transition system
$A$ and a zig-zag of maps $X\stackrel{f}\longleftarrow A
\stackrel{g}\longrightarrow Y$ such that $f$ and $g$ are
$\mathcal{P}$-open.  \ed

That $X$ and $Y$ are $\mathcal{P}$-bisimilar means
that every calculation path $P$ of $\mathcal{P}$ of $X$ is simulated
by a calculation path of $Y$ and vice versa.

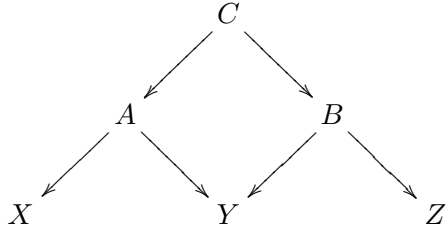
\begin{figure}
\[
\xymatrix{
&& C \ar@{->}[dl] \ar@{->}[dr]&& \\
&A\ar@{->}[dl] \ar@{->}[dr] &&B\ar@{->}[dl] \ar@{->}[dr]& \\
X&&Y&&Z
}
\]
\caption{Bisimulation as an equivalence relation}
\label{equi}
\end{figure}

Bisimilarity is an equivalence relation: it is clearly symmetric, it
is reflexible with $X = A = Y$ and it is transitive since a pullback
of a map satisfying the RLP with respect to a given map still
satisfies the RLP and because of the diagram cartesian in $C$ of
Figure~\ref{equi}.

The following theorem explains the connexion with more usual
($1$-dimensional) notions of bisimulations \cite{MR1365754}.

\bp Take $\mathcal{P} = \{C_1[x]\mid x \in \Sigma\}$. Let $X = (S_X,
\mu:L_X \rightarrow \Sigma,T_X)$ and $Y = (S_Y, \mu:L_Y \rightarrow
\Sigma,T_Y)$ be two cubical transition systems. Then $X$ and $Y$ are
$\mathcal{P}$-bisimilar if and only if there exists a binary relation
$\mathcal{R} \subset S_X \p S_Y$ satisfying the following property:
\begin{enumerate}
\item for every pair $(\alpha,\beta) \in \mathcal{R}$ and every map $c
: C_1[x] \rightarrow X$ with $c(0_1) = \alpha$, there exists a map $d
: C_1[x] \rightarrow Y$ with $d(0_1) = \beta$ and $(c(1_1),d(1_1))
\in \mathcal{R}$
\item for every pair $(\alpha,\beta) \in \mathcal{R}$ and every map $d
: C_1[x] \rightarrow Y$ with $d(0_1) = \beta$, there exists a map $c :
C_1[x] \rightarrow X$ with $c(0_1) = \alpha$ and $(c(1_1),d(1_1)) \in
\mathcal{R}$.
\end{enumerate} \ep

\bpf If $X\stackrel{f}\longleftarrow A \stackrel{g}\longrightarrow Y$
is a map as above, then $\mathcal{R} = \{(f(\alpha),g(\alpha)) \mid
\alpha \hbox{ state of }A\}$ satisfies the two properties of the
statement of the theorem. Conversely, suppose that such a binary
relation $\mathcal{R}$ exists. Let $X\p_\mathcal{R} Y$ be the weak
HDTS with set of states $\mathcal{R}$, with set of actions the one of
$X\p Y$ and such that a transition $(\alpha,u_1,\dots,u_n,\beta)$ of
$X\p Y$ is a transition of $X\p_\mathcal{R} Y$ if and only if $\alpha$
and $\beta$ belong to $\mathcal{R}$. Then consider the image $A$ of
$X\p_\mathcal{R} Y$ by the right adjoint to the inclusion functor
$\cts \subset \whdts$:
\[ A = \liminj_{\begin{array}{c}f= C_n[x_1,\dots,x_n] \rightarrow
X\p_\mathcal{R} Y\\\hbox{ or }f=\dd{x} \rightarrow X\p_\mathcal{R}
Y\end{array}} \dom(f)
\]
Then the composite maps $A\rightarrow X\p_\mathcal{R} Y \rightarrow
X\p Y \rightarrow X$ and $A\rightarrow X\p_\mathcal{R} Y \rightarrow
X\p Y \rightarrow Y$ satisfy the RLP with respect to any map of the
form $\{0_1\} \subset C_1[x]$ for $x \in \Sigma$.  \epf

\bth \label{acc} The class of $\mathcal{P}$-open maps is accessible
and finitely accessibly embedded in the category of maps of cubical
transition systems.  \eth

Note that the following proof does not use the fact that a path is a
cube. It only needs the fact that we consider a \emph{set} of
paths. So our very restrictive choice for the definition of a path
does not matter.

\bpf That it is finitely accessibly embedded (i.e. the inclusion
functor in the category of maps preserves finitely filtered colimits)
comes from the finiteness of the set of states and of the set of
actions of a cube. This class of maps is accessible by
\cite[Proposition~3.3]{MR2506258}.  \epf

Note that the arity of all relation symbols of the theory axiomatizing
the class of $\mathcal{P}$-open maps is finite. This provides another
proof of the fact that the category of $\mathcal{P}$-open maps is
finitely accessibly-embedded (e.g, cf. the proof of
\cite[Theorem~5.9]{MR95j:18001}).

\bth The Bousfield localization of $\bl_{\mathcal{S}}(\cts)$ with
respect to the proper class of $\mathcal{P}$-open maps exists and is a
combinatorial left proper model category. \eth

\bpf The argument is standard. By \cite[Proposition~7.3]{MR1870516},
there exists a regular cardinal $\lambda_1$ such that
$\lambda_1$-filtered colimits of weak equivalences of
$\bl_{\mathcal{S}}(\cts)$ are again weak equivalences. Let $\lambda_2$
be a regular cardinal such that the category of $\mathcal{P}$-open
maps is $\lambda_2$-accessible. Let $\lambda$ be a regular cardinal
sharply bigger than $\lambda_1$ and $\lambda_2$. Consider the
Bousfield localization $\bl_\lambda\bl_{\mathcal{S}}(\cts)$ of
$\bl_{\mathcal{S}}(\cts)$ by a set $\mathcal{A}_\lambda$ of
representatives of the class of $\lambda$-presentable
$\mathcal{P}$-open maps. Then the localization functor
$\bl_\lambda(-)$ is $\lambda$-accessible. Any $\mathcal{P}$-open map
$f$ is a $\lambda$-filtered colimits of maps of $\mathcal{A}_\lambda$,
$f = \liminj_i f_i$ by Theorem~\ref{acc}. So $\bl_\lambda(f) =
\liminj_i \bl_\lambda(f_i)$.  But for every $i$, the map $
\bl_\lambda(f_i)$ is a weak equivalence of
$\bl_{\mathcal{S}}(\cts)$. Therefore $\bl_\lambda(f)$ is a weak
equivalence of $\bl_{\mathcal{S}}(\cts)$ as well. Hence every
$\mathcal{P}$-open map is a weak equivalence of
$\bl_\lambda\bl_{\mathcal{S}}(\cts)$, and therefore the latter model
category is the Bousfield localization. \epf

Note that all maps of $\mathcal{S}$ are actually
$\mathcal{P}$-open. In this new Bousfield localization, two bisimilar
cubical transition systems are weakly equivalent. This new model
category will be the subject of future works. 

Let us conclude this section by mentioning \cite{definissable}. The
class of $\mathcal{P}$-open maps is axiomatized by a set of formulas
such that all quantifiers are bounded. So the latter paper provides
another argument for the existence of the Bousfield localization.

\appendix

\section{Small weak factorization system and coreflectivity}
\label{restriction}

We want to prove in this section that the restriction of a small weak
factorization system to a coreflective locally presentable subcategory
is still small (Theorem~\ref{catlem}) with some additional hypotheses
on the subcategory.

\begin{lem} \label{inc} Let $\LL$ be a coreflective subcategory of a
  cocomplete category $\K$.  Let $I$ be a set of maps of $\K$. One has
  the equality $\inj_\K(I) \cap \Mor(\LL) = \inj_\LL(I)$ and the
  inclusions \[\cell_\LL(I) \subset \cell_\K(I) \cap \Mor(\LL) \subset
  \cof_\K(I) \cap\Mor(\LL) \subset \cof_\LL(I).\] Moreover if $I$ is a
  set of maps of $\LL \subset \K$, then $\cell_\LL(I) = \cell_\K(I)
  \cap \Mor(\LL)$.
\end{lem}

\bpf obvious.  \epf

\begin{lem} \label{l1} (Compare with \cite[Lemma~1.8]{MR1780498}) Let
  $\LL$ be a coreflective subcategory of a locally presentable
  category $\K$.  Let $I$ be a set of maps of $\K$. Let $J$ be a
  solution set for $I$, i.e.  a set of maps of $\LL$ such that every
  map $i \rightarrow w$ of $\Mor(\K)$ from $i \in I$ to $w\in
  \Mor(\LL)$ factors as a composite $i \rightarrow j \rightarrow w$
  with $j\in J$. Then every map $f : X \rightarrow Y$ of $\LL$ can be
  factored as a composite $X \stackrel{g} \longrightarrow P
  \stackrel{h} \longrightarrow Y$ with $g \in \cell_\LL(J)$ and $h\in
  \inj_\LL(I)$.
\end{lem}

\bpf We want to build by transfinite induction on the ordinal $\lambda
\geq 0$ a diagram
\[X =: P_0 \longrightarrow P_1 \longrightarrow \dots \longrightarrow
P_\alpha \longrightarrow P_{\alpha+1} \longrightarrow \dots
\longrightarrow P_\lambda \stackrel{h_\lambda} \longrightarrow Y\]
such that the diagram $P_0 \rightarrow \dots \rightarrow P_\lambda$ is
a transfinite composition of maps belonging to $\cell_\LL(J)$. Since
$P_\lambda$ belongs to $\LL$ and since the category $\LL$ is a full
coreflective subcategory of $\K$, the map $h_\lambda : P_\lambda
\rightarrow Y$ is a map of $\LL$ as well.

Let $P_0 = X$ and $h_0 = f$. For a limit ordinal $\lambda$, let
$P_\lambda = \liminj_{\alpha< \lambda} P_\alpha$. Since the inclusion
functor $\LL \subset \K$ is colimit-preserving, $P_\lambda$ is an
object of $\mathcal{A}$. Let $\lambda \geq 0$ be an ordinal and let us
suppose $P_\alpha$ constructed for $\alpha \leq \lambda$. We want now
to build $P_{\lambda+1}$. Let us consider the set $S_\lambda$ of all
commutative squares
\[
\xymatrix{
  A \fr{} \fd{i} & P_\lambda \fd{h_\lambda} \\
  B \fr{} & Y}
\]
with $i \in I$. The ``density hypothesis'' on $J$ means the existence
of a commutative diagram
\[
\xymatrix{
  A \fr{} \fd{i} & A_s \fr{t_s} \fd{j_s} & P_\lambda \fd{h_\lambda} \\
  B \fr{} & B_s \fr{} & Y}
\]
with $j_s \in J$ (so $A_s$ and $B_s$ both belong to $\LL$), for each
square $s \in S_\lambda$. Let $P_{\lambda+1}$ be the pushout diagram
(in $\LL$ or in $\K$)
\[
\xymatrix{
  \bigsqcup A_s \fD{\bigsqcup j_s} \fR{\bigsqcup_{\{t_s\mid s\in S_\lambda\}} t_s} 
  && P_\lambda \fD{h_{\lambda+1}} \\
  && \\
  \bigsqcup B_s \fR{} && \cocartesien P_{\lambda+1}}
\]
The universal property of the pushout yields a map $h_{\lambda+1} :
P_{\lambda+1} \rightarrow Y$.

Let now $\kappa$ be a regular cardinal exceeding the rank of
presentability of all the objects that occur as domains of maps in
$I$. The required factorization is $X \rightarrow P_\kappa \rightarrow
Y$. Indeed, consider a commutative square of solid arrows of the form
\[
\xymatrix{
A \fr{a} \fd{i} & P_\kappa\fd{} \\
B \fr{}\ar@{-->}[ru]^-{k} & Y}
\]
with $i\in I$. Since $\kappa$ is regular, the diagram $X = P_0
\rightarrow \dots \rightarrow P_\kappa$ is $\kappa$-filtered and since
$\K(A,-)$ commutes with $\kappa$-filtered colimits by hypothesis, the
map $a$ factors as a composite $A \rightarrow P_\lambda \rightarrow
P_\kappa$ for some $\lambda < \kappa$. Let $s\in S_\lambda$ be the
commutative square
\[
\xymatrix{
  A \fr{} \fD{i} & P_\lambda \fd{} \\
& P_\kappa \fd{h_\kappa} \\
  B \fr{} & Y.}
\]
Then the lift $k$ is the bottom composite

\[
\xymatrix{
  A \fr{} \fd{i} & A_s \fr{} \fd{j_s}&\bigsqcup A_s \fd{\bigsqcup j_s} \fr{} & P_\lambda \fd{h_\lambda} & \\
  B \fr{} & B_s \fr{}&\bigsqcup B_s  \fr{} & P_{\lambda+1} \cocartesien \fr{} & P_\kappa.}
\]
\epf 

\begin{lem} \label{l2} Let $\LL$ be a coreflective subcategory of a
  locally presentable category $\K$.  Let $I$ be a set of maps of
  $\K$. Let $J$ be a solution set for $I$ which satisfies $J \subset
  \cof_\K(I)$.  Then there is the equality $\cof_\LL(J) = \cof_\K(I)
  \cap \Mor(\LL)$.
\end{lem}

\bpf One has $\cell_\LL(J) \subset \cof_\K(I)$ since $J\subset
\cof_\K(I)$ and since $\LL$ is coreflective. Since $J$ is a set, every
map of $\cof_\LL(J)$ is a retract of a map of $\cell_\LL(J)$,
therefore $\cof_\LL(J) \subset \cof_\K(I) \cap \Mor(\LL)$. Conversely,
let $f \in \cof_\K(I) \cap \Mor(\LL)$. By Lemma~\ref{l1}, $f$ factors
as a composite
\[
\xymatrix{
\bullet \fD{f}\fR{g} && \bullet \fD{h}\\
&&\\
\bullet \ar@{=}[rr] \ar@{-->}[rruu]^-{k} && \bullet }
\]
with $g\in \cell_\LL(J)$ and $h\in \inj_\LL(I)$. The lift $k$ exists
since $f\in \cof_\K(I)$. The commutative diagram
\[
\xymatrix{
\bullet \ar@{=}[rr]\fD{f} && \bullet \ar@{=}[rr] \fD{g} && \bullet \fD{f}\\
&&&&\\
\bullet \ar@{-->}[rr]^-{k} && \bullet \fR{h} && \bullet}
\]
proves that $f$ is a retract of $g\in \cell_\LL(J)$. Therefore $f\in
\cof_\LL(J)$. Hence the inclusion $\cof_\K(I) \cap \Mor(\LL) \subset
\cof_\LL(J)$.  \epf

We want now conclude the section by giving a sufficient condition for
a small weak factorization system to restrict to a small one on a full
coreflective subcategory. First we recall a definition: 

\bd \cite[Definition~4.14]{MR95j:18001} Let $\K$ be a locally
presentable category. An object $K$ is {\rm injective} with respect to
a cone of maps $(A\rightarrow A_i)_{i\in I}$ if the map $K\rightarrow
\mathbf{1}$ belongs to $\bigcup_{i\in I} \inj(A\rightarrow A_i)$. A
{\rm small cone-injectivity class} is the full subcategory of $\K$ of
objects injective with respect to a given set of cones.  \ed

Hence the conclusion of the section: 

\bth \label{catlem} Let $I$ be a set of maps of a locally presentable
category $\K$.  Let $\LL$ be a coreflective small cone-injectivity
class of $\K$ such that each map of each cone is an element of
$\cof_\K(I)$.  Then there exists a set of maps $J$ of $\mathcal{A}$
such that $\cof_\K(I) \cap \Mor(\LL) = \cof_\LL(J)$.  \eth

\bpf By Lemma~\ref{l2}, it suffices to prove that there exists a set
of maps $J$ of $\mathcal{A}$ which is a solution set for $I$ with $J
\subset \cof_\K(I)$. We mimick the proof of
\cite[Lemma~1.9]{MR1780498}.  Since $\LL$ is a small cone-injectivity
class, it is accessible (and accessibly embedded) by
\cite[Proposition~4.16]{MR95j:18001}. Therefore $\LL$ is locally
presentable by Proposition~\ref{facile_1}.  The inclusion functor
$\Mor(\LL) \subset \Mor(\K)$ is colimit-preserving between two locally
presentable categories (by
\cite[Theorem~2.43]{MR95j:18001}). Therefore it is accessible. So it
satisfies the solution set condition by
\cite[Corollary~2.45]{MR95j:18001}. This means that there exists for
each $i\in I$ a solution set $W_i \subset \Mor(\LL)$, i.e. every map
$i \rightarrow w$ of $\Mor(\K)$ from $i \in I$ to $w\in \Mor(\LL)$
factors as a composite $i \rightarrow w_i \rightarrow w$ for some $w_i
\in W_i$. Consider the set of commutative squares $i \rightarrow w_i$
for $i$ running over the set $I$ and $w_i$ running over the set $W_i$:
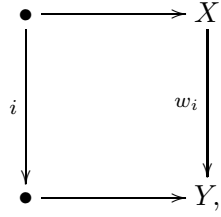
\[
\xymatrix{
\bullet \fD{i} \fR{} && X \fD{w_i} \\
&&\\
\bullet \fR{} && Y,}
\]
Form the pushout diagram
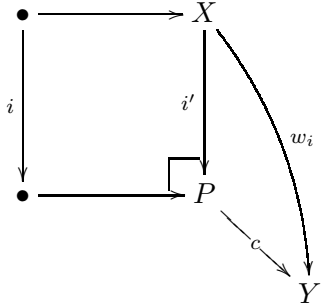
\[
\xymatrix{
\bullet \fD{i} \fR{} && X \fD{i'} \ar@/^10pt/[rddd]^-{w_i}&\\
&&\\
\bullet \fR{} && P \cocartesien \ar@{->}[rd]|-{c}&\\ 
&&& Y}
\]
and factor $c$ as $P \stackrel{p}\rightarrow Q \stackrel{q}
\rightarrow Y$ with $p\in \cell_\K(I)$ and $q \in \inj_\K(I)$. As in
\cite[Lemma~1.9]{MR1780498}, let $J$ be the set of maps $j = pi'$. By
hypothesis, $X$ and $Y$ are cone-injective. Consider a map
$A\longrightarrow Q$ where $A$ is the top of a cone characterizing
$\LL$ as a small cone-injectivity class. Let us consider the
composition
\[A\longrightarrow Q \stackrel{q} \longrightarrow Y.\] Since $Y$ is
cone-injective, there exists a map $A\rightarrow B$ of the cone with
top $A$ and a commutative square of solid arrows of the form
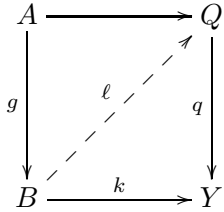
\[
\xymatrix{
  A \fR{} \ar@{->}[dd]_-{g} && Q \fD{q} \\
  && \\
  B \ar@{->}[rr]^-{k} \ar@{-->}[rruu]^-{\ell}&& Y}
\]
Since $g\in \cof_\K(I)$ by hypothesis, and since $q \in \inj_\K(I)$,
the lift $\ell$ exists. This means that $Q$ is cone-injective as well,
i.e. $Q\in \LL$.  Since $\LL$ is a full subcategory of $\K$, we deduce
that $j$ is a map of $\LL$. Therefore, $J\subset \cell_\K(I) \cap
\Mor(\LL)$. Finally, every map $i \rightarrow w$ of $\Mor(\K)$ from $i
\in I$ to $w\in \Mor(\LL)$ factors as a composite $i \rightarrow j
\rightarrow w$ with $j\in J$ by:
\[
\xymatrix{
\bullet\fD{i} \fR{}&& X \fD{j (=pi')}\ar@{=}[rr] && X \fR{} \fD{w_i} &&\bullet \fD{w}  \\
&&&&\\
\bullet \fR{} && Q \fR{q} && Y \fR{}&& \bullet}
\]
\epf

\end{document}


\begin{abstract} Counterexamples for Proposition~8.1 and
  Proposition~8.2 are given. They are used in the paper only to prove
  Corollary~8.3. A proof of this corollary is given without them. The
  proof of the fibrancy of some cubical transition systems is fixed.
\end{abstract}

\maketitle


\begin{propN} (Counterexample for Proposition~8.1) The canonical map
  \[\psi_X:\CSA_1(\cub(X)) \to \cub(\CSA_1(X))\] is bijective on states,
  one-to-one on actions and one-to-one on transitions for all cubical
  transition systems $X$.  There exists a cubical transition system
  $Z$ such that $\psi_Z$ is not surjective on actions and on
  transitions.
\end{propN}

\bpf The map $\psi_X:\CSA_1(\cub(X)) \to \cub(\CSA_1(X))$ is bijective
on states and one-to-one on actions: see the proof of
\cite[Proposition~8.1]{cubicalhdts}. Therefore, it is one-to-one on
transitions by a standard argument already used several times in this
series of papers (see also
\cite[Proposition~4.4]{homotopyprecubical}): if
$(\alpha,u_1,\dots,v_n,\beta)$ and $(\alpha,u'_1,\dots,v'_{n'},\beta)$
are two transitions of $\CSA_1(\cub(X))$ such that
\[(\psi_X(\alpha),\psi_X(u_1),\dots,\psi_X(v_n),\psi_X(\beta)) =
(\psi_X(\alpha),\psi_X(u'_1),\dots,\psi_X(v'_{n'}),\psi_X(\beta)),\]
then $n=n'$ and since $\psi_X$ is one-to-one on states and actions,
one has
\[(\alpha,u_1,\dots,v_n,\beta)=(\alpha,u'_1,\dots,v'_{n'},\beta).\] We
are now going to find a cubical transition system $Z$ such that
$\psi_Z$ is not surjective on actions and on transitions. Consider the
quotient set
\[S=\left(\{\alpha,\beta,\alpha_1,\beta_1,\alpha_2,\beta_2,\alpha_3,\beta_3\}
  \cup \{0,1\}^3\p\{-,+\}\right)/((0,0,0,\pm)=\alpha\hbox{ and
}(1,1,1,\pm)=\beta),\] i.e. with the identifications
$(0,0,0,-)=(0,0,0,+)=\alpha$ and $(1,1,1,-)=(1,1,1,+)=\beta$. Let
$L=\{u_1,u_2,u_3,u'_1,u'_2,u'_3\}$ be a set of actions with the
labelling map defined by $\mu(u_i)=\mu(u'_i)=x_i$ for $i=1,2,3$.  The
cubical transition system $Z$ is intuitively the smallest one having
the set of states $S$ and the set of actions $L$ such that there are
the following maps of cubical transition systems:
\begin{enumerate}
\item The map $C_3[x_1,x_2,x_3]\to Z$ taking the states $(i,j)$ to
  $(i,j,-)$ for all $(i,j)\in\{0,1\}^3$ and the actions $x_i$ to $u_i$
  for $i=1,2,3$ (the actions of $C_3[x_1,x_2,x_3]$ are not denoted by
  $(x_1,1),\dots,(x_3,3)$ because it is already understood that
  $x_1,x_2,x_3$ are distinct).
\item The map $C_3[x_1,x_2,x_3]\to Z$ taking the states $(i,j)$ to
  $(i,j,+)$ for all $(i,j)\in\{0,1\}^3$ and the actions $x_i$ to
  $u'_i$ for $i=1,2,3$.
\item The maps $\cyl(C_1[x_i])\to Z$ for $i=1,2,3$ taking the initial
  state to $\alpha_i$, the final state to $\beta_i$, the action
  $(x_i,0)$ to $u_i$ and the action $(x_i,1)$ to $u'_i$.
\end{enumerate}
It is defined rigorously as the final lift of a cone of maps like as
follows, which always exists since the functor $\omega:\whdts \to
\set^{\{s\}\cup \Sigma}$ forgetting the transitions is topological:
\begin{multline*}
\omega(C_3[x_1,x_2,x_3])\to
(\{0,1\}^3\p\{-\},\{u_1\},\{u_2\},\{u_3\}) \subset W \\
\omega(C_3[x_1,x_2,x_3])\to
(\{0,1\}^3\p\{+\},\{u'_1\},\{u'_2\},\{u'_3\}) \subset W\\\omega
(\cyl(C_1[x_1]))\to (\{\alpha_1,\beta_1\},\{u_1,u'_1\},\varnothing,\varnothing) \subset W\\
\omega(\cyl(C_1[x_2]))\to (\{\alpha_2,\beta_2\},\varnothing,\{u_2,u'_2\},\varnothing) \subset W\\
\omega(\cyl(C_1[x_3]))\to (\{\alpha_3,\beta_3\},\varnothing,\varnothing,\{u_3,u'_3\}) \subset W
\end{multline*}
with $W= (S,\{u_1,u'_1\},\{u_2,u'_2\},\{u_3,u'_3\})$. One has
$\omega(Z)=W$. The weak HDTS $Z$ is cubical since $\cts$ is a
coreflective subcategory of the category of weak HDTS.
The key fact is that $Z$ contains the transitions
$(\alpha_i,u_i,\beta_i)$ and $(\alpha_i,u'_i,\beta_i)$ for
$i=1,2,3$. Therefore the canonical map $\phi_Z:Z\to\CSA_1(Z)$
identifies the actions $u_i$ and $u'_i$ for $i=1,2,3$:
$\phi_Z(u_i)=\phi_Z(u'_i)=x_i$ for $i=1,2,3$.  So $\CSA_1(Z)$ contains
the five transitions (remember that $\phi_Z$ is bijective on states)
\[
(\alpha,x_1,x_2,x_3,\beta),(\alpha,x_1,(1,0,-)),((1,0,-),x_2,x_3,\beta),(\alpha,x_1,x_2,(1,1,+)),((1,1,+),x_3,\beta).
\]
By the composition axiom, $\CSA_1(Z)$ contains the transition
$((1,0,-),x_2,(1,1,+))$ which corresponds to a unique map $C_1[x_2]\to
\CSA_1(Z)$.  Hence the cubical transition system $\cub(\CSA_1(Z))$
contains a transition from $(1,0,-)$ to $(1,1,+)$ indexed by an action
$u''_2$ labelled by $x_2$ which is distinct from $u_2$ and $u'_2$.
The point is that in $\cub(Z)$, the transition
$(\alpha_i,u_i,\beta_i)$ becomes a transition $(\alpha_i,v_i,\beta_i)$
and the transition $(\alpha_i,u'_i,\beta_i)$ becomes a transition
$(\alpha_i,v'_i,\beta_i)$ with
$\mu(u_i)=\mu(v_i)=\mu(u'_i)=\mu(v'_i)=x_i$ for $i=1,2,3$.  So the
canonical map $\phi_{\cub(Z)}:\cub(Z)\to \CSA_1(\cub(Z))$ does not
identify the actions $u_i$ and $u'_i$ for $i=1,2,3$. Therefore the
composition axiom cannot be applied in $\CSA_1(\cub(Z))$ to create a
transition from $(1,0,-)$ to $(1,1,+)$. Hence the map
$\psi_Z:\CSA_1(\cub(Z)) \to \cub(\CSA_1(Z))$ is not surjective on
actions and on transitions. \epf

\begin{propN} (Counterexample for Proposition~8.2) There exists a
  weak equivalence of $\cts$ (the left determined model structure)
  such that $\cub(f)$ is not a weak equivalence of $\cts$. \end{propN}

\bpf Let $f=\phi_Z:Z \to \CSA_1(Z)$ with $Z$ as above. Then $f$ is a
weak equivalence of $\cts$ by \cite[Theorem~7.10]{cubicalhdts} since
$\CSA_1(f):\CSA_1(Z)\to \CSA_1(\CSA_1(Z))$ is an isomorphism. The
source of $\CSA_1(\cub(f))$ is $\CSA_1(\cub(Z))$. The target of
$\CSA_1(\cub(f))$ is $\CSA_1(\cub(\CSA_1(Z)))$ which is equal to
$\cub(\CSA_1(Z))$ because the latter satisfies CSA1 (see the beginning
of the proof of \cite[Proposition~8.1]{cubicalhdts}). Therefore
$\CSA_1(\cub(f))$ cannot be an isomorphism and $\cub(f)$ is not a weak
equivalence of $\cts$ by \cite[Theorem~7.10]{cubicalhdts}.  \epf

\begin{propN} (Corollary~8.3 fixed) Every weak equivalence of $\cts$
  belongs to $\mathcal{W}_{\cub}$.
\end{propN}

\bpf The class of maps $\mathcal{W}_{\cub}$ is, by definition, the
localizer generated by the maps of cubical transition systems $f$ such
that $\cub(f)$ is a weak equivalence of $\cts$. This localizer
contains the smallest one, which is precisely the class of weak
equivalences of the left determined model structure $\cts$.  \epf

Let $\I$ be the set of generating cofibrations of $\cts$. Let $S$ be
an arbitrary set of maps in $\cts$.  It is claimed in
\cite{cubicalhdts} that the class of fibrant objects of the Bousfield
localization by the set of maps $S$ of the left determined model
structure $\cts$ is the class of $\Lambda_{\cts}(V,S,\I)$-injective
objects. Using Olschok's theorems, it is only possible to say that the
class of fibrant objects is the class of
$\Lambda_{\cts}(V,S^{op},\I)$-injective objects where $S^{cof}$ is a
set of cofibrant replacements for the maps of $S$. Since
$\varnothing=\varnothing^{cof}$, it is correct to say that the class
of fibrant objects of the left determined model structure of $\cts$ is
the class of $\Lambda_{\cts}(V,\varnothing,\I)$-injective objects. So
the proof of Proposition~7.8 is correct. However,
``$\Lambda_{\cts}(V,S,\I)$-injective'' must be replaced by
``$\Lambda_{\cts}(V,S^{op},\I)$-injective'' page 318 before and in the
proof of Theorem~6.3. And the proofs of Theorem~8.10 and Theorem~8.11
must be modified. More precisely, the proof of the following fact must
be modified, and without using Theorem~8.10 (the characterization of
the weak equivalences of $\bl_{\mathcal{S}}(\cts)$) and Theorem~8.11
to avoid any vicious circle:

\begin{propN} (Proof of fibrancy fixed) Let
  $\mathcal{S}=\{p_x:C_1[x]\sqcup C_1[x] \to \dd{x}\mid x\in
  \Sigma\}$.  Then any $\mathcal{S}$-injective cubical transition
  system is fibrant in the Bousfield localization
  $\bl_{\mathcal{S}}(\cts)$ of the left determined model structure of
  $\cts$ by the set of maps $\mathcal{S}$. \end{propN}

\bpf This is \cite[Proposition~8.4]{homotopyprecubical}. \epf
